\documentclass [12pt]{article}
\usepackage{amsmath,amssymb}
\begin{document}

\title{ Cramer rule over quaternion skew field}
\author {Ivan Kyrchei  \footnote{Pidstrygach Institute for Applied Problems of Mechanics and Mathematics of NAS of Ukraine,
str. Naukova 3b, Lviv, Ukraine, 79053, kyrchei@lms.lviv.ua}}

\date{}
 \maketitle


\begin{abstract}
 New definitions of
determinant functionals over the quaternion skew field are given
in this paper. The inverse matrix over the quaternion skew field
is represented by analogues of the classical adjoint matrix.
Cramer rule for  right and left quaternionic systems of linear
equations have been obtained.


\end{abstract}

\section{Introduction}

 For a representation of solution of a system of linear equations over the quaternion
 skew field $\bf H$
 by Cramer rule
  is necessary to represent the  inverse  of the coefficient matrix  by the classical
   adjoint matrix.
  The crucial importance for this has a definition of the  determinant of a square
   matrix over $\bf H$.
On the whole the theory of determinants of matrices with
noncommutative entries, (which are also  defined as noncommutative
determinants),
 can be divided into three   methods. Let ${\rm M}\left( {n,\bf K} \right)$ be
  the ring of $n\times n$ matrices with entries in a
ring $ {\bf K}$. The first approach \cite{As1,co4,dy5} to defining
the determinant of a matrix in ${\rm M}\left( {n,\bf K} \right)$
is as follows.
\newtheorem{definition}{Definition}[section]
 \begin{definition} Let a
functional ${\rm d}:M\left( {n,\bf K} \right) \to \bf K$ satisfy
the following axioms.
\begin{itemize}
\item [] {\bf Axiom 1} ${\rm d}\left( {{\rm {\bf A}}} \right) = 0$ if and only if the matrix ${\rm {\bf A}}$ is singular.
\item [] {\bf Axiom 2} ${\rm d}\left( {{\rm {\bf A}} \cdot
{\rm {\bf B}}} \right) = {\rm d}\left( {{\rm {\bf A}}} \right)
\cdot {\rm d}\left( {{\rm {\bf B}}} \right)$ for $\forall {\rm
{\bf B}}\in {\rm M}\left( {n,\bf K} \right)$.
\item [] {\bf Axiom 3} If the matrix ${\rm {\bf A}}'$ is obtained
from ${\rm {\bf A}}$ by adding a left-multiple of a row to another
row or a right-multiple of a column to another column, then ${\rm
d}\left( {{\rm {\bf A}}}\right)'={\rm d}\left( {{\rm {\bf
A}}}\right)$.
\end{itemize}
Then a value of the functional ${\rm d}$ is called the determinant
of the matrix ${\rm {\bf A}}\in {\rm M}\left( {n,\bf K}
\right)$.\end{definition}
Examples of such determinant are the
 determinants of Study and  Diedonn\'{e}. If a determinant functional satisfies  Axioms 1, 2,
 3, then it takes on a value in a commutative subset of the ring.
 It is proved in
\cite{As1}.  Therefore a determinant representation of an inverse
matrix by such determinants is impossible.   This  reason  compels
to  define determinant functionals unsatisfying   all the above
axioms. However Axiom 1 is considered \cite{dy5} indispensable for
the utility of the notion of a determinant.

In another way   a noncommutative determinant  is defined as a
rational function from  entries of a matrix. Herein I. M. Gelfand
and V. S. Retah  have reached  the greatest success by the theory
of quasideterminants  \cite{ge6,ge7}. An arbitrary $n\times n$
matrix over a skew field is associated with an $n\times n$ matrix
whose entries are quasideterminants. I. M. Gelfand
 and V. S. Retah   transfer from  a commutative case not the
concept of a determinant but its relations to minors. Since
quasideterminants can not be expanded  by cofactors along an
arbitrary row or column, an inverse matrix is not represented by
the adjoint classical matrix in this case as well.

At last, at the third approach  a noncommutative determinant is
defined as the alternating sum of $n!$ products of entries of a
matrix but by specifying a certain ordering of coefficients in
each term.  E. H. Moore was the first who
 achieved the fulfillment of the
main Axiom 1 by such definition of  a noncommutative determinant.
This is done  not for all square matrices over a skew field but
rather only Hermitian matrices. He has defined the determinant of
a Hermitian matrix ${\rm {\bf A}}=(a_{ij})_{n\times n} $, ( i.e.
 $a_{ij}=\overline{a_{ji}}$), over a skew field with an
involution by induction on $n$  in the following way (\cite{dy5}).

Denote by $ {{\rm {\bf A}}(i \to j )} $    the matrix obtained
from ${\rm {\bf A}}$  by replacing its $j$th column with the $i$th
column, and then by deleting both the
 $i$th row and column. By definition, put
\begin{equation}\label{kyr1}
\begin{array}{cc}
  {\rm Mdet}  {\rm{\bf A}} = {\left\{
{{\begin{array}{*{20}c}
 {a_{11} {\rm ,}\quad \quad \quad \quad \quad \quad \quad \quad \,\,\,\,\, \quad n = 1}
\hfill \\
 {{\sum\limits_{j = 1}^{n} {\varepsilon _{ij} a_{ij} {\rm Mdet} \left( {{\rm {\bf
A}}{\rm (}i \to j{\rm )}} \right),\,\,n > 1}} } \hfill \\
\end{array}} } \right.}
\end{array}
\end{equation}
\noindent where $\varepsilon _{kj} = {\left\{
{{\begin{array}{*{20}c}
 {1,\,\,\,i = j} \hfill \\
 { - 1,\,\,i \ne j} \hfill \\
\end{array}} } \right.}$. Another definition of this determinant is represented
(\cite{As1}) in terms of permutations:
\[
{\rm{Mdet}}\, {\rm {\bf A}} = {\sum\limits_{\sigma \in S_{n}}
{{\left| {\sigma} \right|}{a_{n_{11} n_{12}}  \cdot \ldots \cdot
a_{n_{1l_{1}}  n_{11}}  \cdot }} {a_{n_{21} n_{22}}} \cdot \ldots
\cdot {a_{n_{rl_{1}}  n_{r1}}}}.
\]
 The disjoint
cycle representation of the permutation $\sigma \in S_{n}$ is
written in the normal form, $ \sigma = \left( {n_{11} \ldots
n_{1l_{1}} }  \right)\left( {n_{21} \ldots n_{2l_{2}} }
\right)\ldots \left( {n_{r1} \ldots n_{rl_{r}} }  \right). $

However there was no extension of the definition of the Moore
determinant to arbitrary square matrices. F. J. Dyson  has
emphasized this point  in \cite{dy5}. Longxuan Chen has offered
the following decision of this problem in \cite{ch2,ch3}. He has
defined the determinant of an arbitrary square matrix ${\rm {\bf
A}}=(a_{ij}) \in {\rm M}\left( {n,\bf{H}} \right)$ over the
quaternion skew field $\bf{H}$ as follows.
\[
\begin{array}{c}
  \det {\rm {\bf A}} = {\sum\limits_{\sigma \in S_{n}}  {\varepsilon
\left( {\sigma}  \right)a_{n_{1} i_{2}}  \cdot a_{i_{2} i_{3}}
\ldots \cdot a_{i_{s} n_{1}}  \cdot} } \ldots \cdot a_{n_{r}
k_{2}}  \cdot \ldots \cdot a_{k_{l} n_{r}}, \\
   \sigma = \left( {n_{1} i_{2} \ldots i_{s}}  \right)\ldots \left(
{n_{r} k_{2} \ldots k_{l}}  \right),\\
  n_{1} > i_{2} ,i_{3} ,\ldots ,i_{s} ;\ldots ;n_{r} > k_{2} ,k_{3}
,\ldots ,k_{l}, \\
  n = n_{1} > n_{2} > \ldots > n_{r} \ge 1.
\end{array}
\]
 L. Chen  has obtained a determinant
representation of an inverse matrix over the quaternion skew field
even though the determinant does not satisfy Axiom 1.
  However this determinant
also can not be expanded  by cofactors along an arbitrary row or
column with the exception of the $n$th row. Therefore he  has not
obtained the classical adjoint matrix or its analogue as well.

If ${\left\| {{\rm {\bf A}}} \right\|}:=\det({\rm {\bf A}}^{*}{\rm
{\bf A}}) \ne 0$ for ${\rm {\bf A}} = \left( {\alpha _{1} ,\ldots
,\alpha _{m}} \right)$ over {\bf H}, then $\exists {\rm {\bf A}}^{
- 1} = \left( {b_{jk}} \right)$, where
\[
\overline {b_{jk}}  = {\frac{{1}}{{{\left\| {{\rm {\bf A}}}
\right\|}}}}\omega _{kj}, \quad \left( {j,k = 1,2,\ldots ,n}
\right),
\]
\[
\omega _{kj} = \det \left( {\alpha _{1} \ldots \alpha _{j - 1}
\alpha _{n} \alpha _{j + 1} \ldots \alpha _{n - 1} \delta _{k}}
\right)^{ *} \left( {\alpha _{1} \ldots \alpha _{j - 1} \alpha
_{n} \alpha _{j + 1} \ldots \alpha _{n - 1} \alpha _{j}} \right).
\]
\noindent Here $\alpha _{i} $ is the $i$th column of ${\rm {\bf
A}}$, $\delta _{k} $ is the $n$-dimension column with 1 in the
$k$th row and 0 in others. He defined ${\left\| {{\rm {\bf A}}}
\right\|}:=\det({\rm {\bf A}}^{*}{\rm {\bf A}})$ as the double
determinant. If ${\left\| {{\rm {\bf A}}} \right\|} \ne 0$, then a
solution of a right system of linear equations ${\sum\nolimits_{j
= 1}^{n} {\alpha _{j} x_{j} = \beta} } $ over {\bf H}  is
represented by the following formula,  defined as Cramer formula,
\[
x_{j} = {\left\| {{\rm {\bf A}}} \right\|}^{ - 1}\overline {{\rm
{\bf D}}_{j}}, \,\forall j=\overline{1,n},
\]
\noindent where
\[
{\rm {\bf D}}_{j} = \det \left( {{\begin{array}{*{20}c}
 {\alpha _{1}^{ *} }  \hfill \\
 { \vdots}  \hfill \\
 {\alpha _{j - 1}^{ *} }  \hfill \\
 {\alpha _{n}^{ *} }  \hfill \\
 {\alpha _{j + 1}^{ *} }  \hfill \\
 {\vdots} \hfill \\
 {\alpha _{n - 1}^{ *} }  \hfill \\
 {\beta ^{ *} } \hfill \\
\end{array}} } \right)\left( {{\begin{array}{*{20}c}
 {\alpha _{1}}  \hfill & {\ldots}  \hfill & {\alpha _{j - 1}}  \hfill &
{\alpha _{n}}  \hfill & {\alpha _{j + 1}}  \hfill & {\ldots}
\hfill & {\alpha _{n - 1}}  \hfill & {\alpha _{j}}  \hfill \\
\end{array}} } \right).
\]
Here $\alpha _{i} $ is the $i$th column of ${\rm {\bf A}}$,
$\alpha _{i}^{ *}  $ is the $i$th row of ${\rm {\bf A}}^{ *} $,
and $\beta ^{ *} $ is the $n$-dimension row vector conjugated with
 $\beta $.

The row and  column determinants of a square matrix over the
quaternion skew field are defined in this work. (For the first
time we introduced these definitions in \cite{ky8}). Their
properties of an arbitrary square matrix  and Hermitian over the
quaternion skew field are investigated. The determinant
representations of an inverse matrix  by  analogues of the adjoint
matrix are obtained. Generalizations of Cramer rule for left and
right systems of linear equations over the quaternion skew field
are obtained as well.
\section{ Definitions and basic properties\\ of the column and  row determinants}
Throughout  this article the skew field $\bf H $ is the quaternion
division algebra  generated by four basic elements $1, i, j, k$
over the field of real numbers $\bf R$ with the famous Hamilton's
relations $ i^{2}=j^{2}=k^{2}=ijk=-1. $ Define
$q_{n}=w_{n}+x_{n}i+y_{n}j+z_{n}k \in \bf H$ for $n=(1,2)$.
\textit{Addition} and \textit{subtraction} of quaternions is
defined by $ q_{1} \pm q_{2} = (w_{1} \pm w_{2}) + (x_{1} \pm
x_{2})i + (y_{1} \pm y_{2})j + (z_{1} \pm z_{2})k$.
\textit{Multiplication} of quaternions is defined by
\[\begin{array}{c}
   q_{1}q_{2} =   (w_{1}w_{2} - x_{1}x_{2} - y_{1}y_{2} - z_{1}z_{2})
   +(w_{1}x_{2} + x_{1}w_{2} + y_{1}z_{2} - z_{1}y_{2})i+\\
  +(w_{1}y_{2}
- x_{1}z_{2} + y_{1}w_{2} + z_{1}x_{2})j
   +(w_{1}z_{2} + x_{1}y_{2}
- y_{1}x_{2} + z_{1}w_{2})k.
\end{array}
\]
 The \textit{conjugate} of a quaternion $q=w+xi+yj+zk$ is defined to be
$\overline{q}=w-xi-yj-zk$ at that $\overline {p + q} = \overline
{q} + \overline {p} $, $\overline {p \cdot q} = \overline {p}
\cdot \overline {q} $, $\overline {\overline {q}} = q$ for all
$q,p \in \bf H$. The \textit{norm} of a quaternion is defined by
${\rm n}(q)=w^{2}+x^{2}+y^{2}+z^{2}$. The norm is a real-valued
function, and the norm of a product of quaternions satisfies the
properties ${\rm n}\left( {p \cdot q} \right) = {\rm n}\left( {p}
\right) \cdot {\rm n}\left( {q} \right)$ and ${\rm n}\left(
\overline{q} \right) = {\rm n}(q)$. The \textit{trace} of a
quaternion is defined by ${\rm t}\left( {q} \right) = q +
\overline {q}$. The trace is a real-valued function as well. The
trace of a product of quaternions satisfies the rearrangement
property ${\rm t}\left( {q \cdot p} \right) = {\rm t}\left( {p
\cdot q} \right)$.
\begin{definition}Suppose $S_{n}$ is the symmetric group on the set
$I_{n}=\{1,\ldots,n\}$. We say that the permutation $\sigma \in
S_{n}$ is written by the direct product of disjoint cycles if its
subscription by usual two-line representation corresponds to the
cycle notation, i.e.
\begin{equation}\label{kyr2}
\sigma =
\begin{pmatrix}
  n_{11} & n_{12} & \ldots & n_{1l_{1}} & \ldots & n_{r1} & n_{r2} & \ldots & n_{rl_{r}}\\
  n_{12} & n_{13} & \ldots & n_{11} & \ldots & n_{r2} & n_{r3} & \ldots & n_{r1}
\end{pmatrix}
\end{equation}
\end{definition}
\begin{definition}
We say that the cycle notation of the permutation $\sigma \in
S_{n}$ is left ordered if the first elements from the left in each
cycles are elements  bringing  closure to the cycles. This means
that if the subscription of the permutation $\sigma \in S_{n}$ by
the direct product of disjoint cycles has the form (\ref{kyr2}),
then the left-ordered cycle notation is represented by
\[
\sigma = \left( {n_{11} n_{12} \ldots n_{1l_{1}} } \right)\,\left(
{n_{21} n_{22} \ldots n_{2l_{2}} }  \right)\ldots \left( {n_{r1}
n_{r2} \ldots n_{rl_{r}} }  \right)
\]
\end{definition}
\begin{definition}
We say that the cycle notation of the permutation $\sigma \in
S_{n}$ is right ordered if the first elements from the right in
each cycles are elements  bringing  closure to the cycles. This
means that if the subscription of the permutation $\sigma \in
S_{n}$ by the direct product of disjoint cycles has the form
(\ref{kyr2}), then the right-ordered cycle notation is represented
by
\[
\sigma = \left(  n_{12} \ldots n_{1l_{1}}{n_{11} } \right)\,\left(
 n_{22} \ldots n_{2l_{2}}{n_{21} }  \right)\ldots \left(
n_{r2} \ldots n_{rl_{r}} {n_{r1}}  \right)
\]
\end{definition}
\begin{definition}\label{kyrc1}
The $i$th row determinant of  ${\rm {\bf A}}=(a_{ij})\in {\rm
M}\left( {n,\bf H} \right)$
 is defined
as the alternating sum of  $n!$ products  of entries of ${\rm {\bf
A}}$, during which the index permutation of every product is
written by the direct product of disjoint cycles. If the
permutation is even, then  product  of entries has a sign "$+$".
If the permutation is odd, then  product  of entries has a sign
"$-$". That is
\[
{\rm{rdet}}_{ i} {\rm {\bf A}} = {\sum\limits_{\sigma \in S_{n}}
{\left( { - 1} \right)^{n - r}{a_{i{\kern 1pt} i_{k_{1}}} }
{a_{i_{k_{1}}   i_{k_{1} + 1}}}   \ldots } } {a_{i_{k_{1} + l_{1}}
 i}}  \ldots  {a_{i_{k_{r}}  i_{k_{r} + 1}}}
\ldots  {a_{i_{k_{r} + l_{r}}  i_{k_{r}} }},
\]
where  $S_{n}$ is the symmetric group on the set $I_{n}$. The
left-ordered cycle notation of the permutation $\sigma$ is written
as follows
\[
   \sigma = \left( {i\,i_{k_{1}}  i_{k_{1} + 1} \ldots i_{k_{1} +
l_{1}} } \right)\left( {i_{k_{2}}  i_{k_{2} + 1} \ldots i_{k_{2} +
l_{2}} } \right)\ldots \left( {i_{k_{r}}  i_{k_{r} + 1} \ldots
i_{k_{r} + l_{r}} } \right).
\]
Here the index $i$ opens the first cycle from the left  and other
cycles satisfy the following conditions
 \[
 i_{k_{2}}  < i_{k_{3}}  < \ldots < i_{k_{r}},
\quad i_{k_{t}}  < i_{k_{t} + s}, \quad \left( {\forall t =
\overline {2,r}}  \right), \quad \left( {\forall s = \overline
{1,l_{t}} }  \right).
\]
\end{definition}

We shall further consider  ${\rm{rdet}}_{i} {\rm {\bf A}}
 \left( {\forall i = \overline {1,n}}  \right)$ as a sum of $n!$
 monomials whose coefficients are entries of ${\rm {\bf A}}$.

  Let ${\rm
{\bf a}}_{.j} $ be the $j$th column and ${\rm {\bf a}}_{i.} $ be
the $i$th row of a matrix ${\rm {\bf A}}\in {\rm M}\left( {n,\bf
H} \right)$. Denote by ${\rm {\bf A}}_{.j} \left( {{\rm {\bf b}}}
\right)$  the matrix obtained from ${\rm {\bf A}}$ by replacing
its $j$th column with the column ${\rm {\bf b}}$ , and by ${\rm
{\bf A}}_{i.} \left( {{\rm {\bf b}}} \right)$ denote the matrix
obtained from ${\rm {\bf A}}$ by replacing its $i$th row with the
row ${\rm {\bf b}}$. Denote by ${\rm {\bf A}}_{}^{i{\kern 1pt} j}
$
 the submatrix of ${\rm {\bf A}}$ obtained by deleting both
the $i$th row and the $j$th column. The following lemma enables us
to expand
 ${\rm{rdet}}_{{i}}\, {\rm {\bf A}}$ $(\forall {i = \overline {1,n}} )$
 by cofactors along  the $i$-th row. The calculation of the row
determinant of a $n\times n$ matrix is reduced to the calculation
of the row determinant of a lower dimension matrix.

\newtheorem{lemma}{Lemma}[section]
\begin{lemma} \label{kyrc2} Let $R_{i{\kern 1pt} j}$ be the right $ij$-th
cofactor of ${\rm {\bf A}}\in {\rm M}\left( {n,\bf H} \right)$,
that is $ {\rm{rdet}}_{{i}}\, {\rm {\bf A}} = {\sum\limits_{j =
1}^{n} {{a_{i{\kern 1pt} j} \cdot R_{i{\kern 1pt} j} } }} $, $(
\forall {i = \overline {1,n}})$.  Then
\[
 R_{i{\kern 1pt} j} = {\left\{ {{\begin{array}{*{20}c}
  - {\rm{rdet}}_{{j}}\, {\rm {\bf A}}_{{.{\kern 1pt} j}}^{{i{\kern 1pt} i}} \left( {{\rm
{\bf a}}_{{.{\kern 1pt} {\kern 1pt} i}}}  \right),& {i \ne j},
\hfill \\
 {\rm{rdet}} _{{k}}\, {\rm {\bf A}}^{{i{\kern 1pt} i}},&{i = j},
\hfill \\
\end{array}} } \right.}
\]
 where  ${\rm {\bf A}}_{.{\kern 1pt} j}^{i{\kern 1pt} i} \left(
{{\rm {\bf a}}_{.{\kern 1pt} {\kern 1pt} i}}  \right)$ is obtained
from ${\rm {\bf A}}$   by replacing the $j$th column with the
$i$th column, and then by deleting both the $i$th row and column;
$k = \min {\left\{ {I_{n}}  \right.} \setminus {\left. {\{i\}}
\right\}} $.
\end{lemma}
{\textit{Proof.}} First  we  prove that $R_{i{\kern 1pt} i} =
{\rm{rdet}} _{{k}}\, {\rm {\bf A}}^{{i{\kern 1pt} i}}$, $k = \min
{\left\{ {I_{n}}  \right.} \setminus {\left. {\{i\}} \right\}} $.
If $i = 1$, then ${\rm{rdet}} _{1}\, {\rm {\bf A}} = {a_{11} \cdot
R{}_{11} + a_{12} \cdot R{}_{12} + \ldots + a_{1n} \cdot
R{}_{1n}}$. Consider monomials of ${\rm{rdet}} _{1}\, {\rm {\bf
A}}$ such that the coefficient $a_{11} $ is the first from the
left in each of  their:
\[
\begin{array}{c}
   a_{11} \cdot R_{11} = {\sum\limits_{\tilde {\sigma}  \in S_{n}}
{\left( { - 1} \right)^{n - r}a_{11}  a_{2{\kern 1pt} i_{k_{2}} }
 \ldots } }a_{i_{k_{2} + l_{2}} {\kern 1pt} 2}  \ldots
 a_{i_{k_{r}} {\kern 1pt} i_{k_{r} + 1}}  \ldots
a_{i_{k_{r} + l_{r}}  {\kern 1pt} i_{k_{r}} },\\
  \tilde {\sigma}  = \left( {1} \right)\left( {2\,i_{k_{2}}  \ldots
i_{k_{2} + l_{2}} }  \right)\ldots \left( {i_{k_{r}}  i_{k_{r} +
1} \ldots i_{k_{r} + l_{r}} }  \right).
\end{array}
\]
By factoring the common left-side factor $a_{11} $, we obtain
\[
\begin{array}{c}
   a_{11}  R_{11}
 = a_{11}  {\sum\limits_{\tilde {\sigma} _{1} \in S_{n-1} } {\left( { - 1} \right)^{n - 1 - \left( {r - 1}
\right)}a_{2{\kern 1pt}  i_{k_{2}} }  \ldots} } a_{i_{k_{2} +
l_{2}}   2} \ldots a_{i_{k_{r}}  i_{k_{r} + 1}} \ldots a_{i_{k_{r} + l_{r} } i_{k_{r}} },\\
  \tilde {\sigma} _{1} = \left( {2\,i_{k_{2}}  \ldots i_{k_{2} + l_{2}} }
\right)\ldots \left( {i_{k_{r}}  i_{k_{r} + 1} \ldots i_{k_{r} +
l_{r}} } \right).
\end{array}
\]
Here  $S_{n-1}$ is the symmetric group on $I_{n}\setminus\{1\}$.
The numbers of the disjoint cycles and  coefficients of every
monomial of $R_{11}$  decrease  by one. An element of the second
row opens  each monomial of $R_{11}$ on the left. There are no
elements of the first row and  column of ${\rm {\bf A}}$ among its
coefficients. Thus, we have
\begin{equation}\label{kyr3}
 R_{11} ={\sum\limits_{\tilde {\sigma} _{1} \in S_{n-1} } {\left( { - 1}
\right)^{n - 1 - \left( {r - 1} \right)}a_{2{\kern 1pt}  i_{k_{2}}
}  \ldots} } a_{i_{k_{2} + l_{2}}   2} \ldots a_{i_{k_{r} + l_{r}
}  i_{k_{r}} }  = {\rm{rdet}} _{2} {\rm {\bf A}}^{11}.
\end{equation}
If now $i \ne 1$, then
\begin{equation}\label{kyr4}
 {\rm{rdet}}_{{i}}\, {\rm {\bf A}} ={ a_{i1} \cdot R_{i1} +
a_{i2} \cdot R{}_{i2} + \ldots + a_{i{\kern 1pt} n} \cdot
R{}_{i{\kern 1pt} n}}
\end{equation}
Consider monomials of ${\rm{rdet}}_{{i}}\, {\rm {\bf A}}$ such
that $a_{i{\kern 1pt} i} $ is the first from the left in each of
their:
\[
\begin{array}{c}
   a_{i{\kern 1pt} i} \cdot R_{i{\kern 1pt} i} =
{\sum\limits_{\mathord{\buildrel{\lower3pt\hbox{$\scriptscriptstyle\frown$}}\over
{\sigma} } \in \,S_{n}}  {\left( { - 1} \right)^{n - r}a_{i{\kern
1pt} i}  a_{1{\kern 1pt}  i_{k_{2}} }   \ldots } }a_{i_{k_{2} +
l_{2}}  1}  \ldots  a_{i_{k_{r}}  i_{k_{r} + 1}} \ldots
a_{i_{k_{r} + l_{r}}  i_{k_{r}} },\\
  \mathord{\buildrel{\lower3pt\hbox{$\scriptscriptstyle\frown$}}\over {\sigma
}} = \left( {i} \right)\left( {1\,i_{k_{2}}  \ldots i_{k_{2} +
l_{2}} } \right)\ldots \left( {i_{k_{r}}  i_{k_{r} + 1} \ldots
i_{k_{r} + l_{r}} } \right)
\end{array}
\]
Again by factoring the common left-side factor $a_{i{\kern
1pt}i}$, we get
\[
\begin{array}{c}
  a_{i{\kern 1pt} i} \cdot R_{i{\kern 1pt} i} = a_{i{\kern 1pt} i}
\cdot
{\sum\limits_{\mathord{\buildrel{\lower3pt\hbox{$\scriptscriptstyle\frown$}}\over
{\sigma} } _{1} \in
\mathord{\buildrel{\lower3pt\hbox{$\scriptscriptstyle\frown$}}\over{S}_{n-1}}}
{\left( { - 1} \right)^{n - 1 - \left( {r - 1} \right)}a_{1 {\kern
1pt} i_{k_{2} }}  \ldots} } a_{i_{k_{2} + l_{2}}   1}  \ldots
a_{i_{k_{r} + l_{r}}   i_{k_{r}} }, \\
  \mathord{\buildrel{\lower3pt\hbox{$\scriptscriptstyle\frown$}}\over {\sigma
}} _{1} = \left( {1\,i_{k_{2}}  \ldots i_{k_{2} + l_{2}} }
\right)\ldots \left( {i_{k_{r}}  i_{k_{r} + 1} \ldots i_{k_{r} +
l_{r}} }  \right).
\end{array}
\]
Here
$\mathord{\buildrel{\lower3pt\hbox{$\scriptscriptstyle\frown$}}\over{S}_{n-1}}$
is the symmetric group on $I_{n}\setminus\{i\}$. The numbers of
disjoint cycles and coefficients of every monomial of $R_{i{\kern
1pt} i}$ again decrease  by  one. An element of the first row
  opens  on the left each
monomial of $R_{i{\kern 1pt} i}$, and there are no elements of the
$i$th row and  column of ${\rm {\bf A}}$ among its coefficients.
Thus, we obtain
\begin{equation}\label{kyr5}
R_{i{\kern 1pt} i}
={\sum\limits_{\mathord{\buildrel{\lower3pt\hbox{$\scriptscriptstyle\frown$}}\over
{\sigma} } _{1} \in
\mathord{\buildrel{\lower3pt\hbox{$\scriptscriptstyle\frown$}}\over{S}_{n-1}}}
{\left( { - 1} \right)^{n - 1 - \left( {r - 1} \right)}a_{1 {\kern
1pt} i_{k_{2} }}  \ldots} } a_{i_{k_{2} + l_{2}}   1} \ldots
a_{i_{k_{r} + l_{r}}   i_{k_{r}} }  = {\rm{rdet}} _{1} {\rm {\bf
A}}^{i{\kern 1pt} i}.
\end{equation}
Combining (\ref{kyr3}) and (\ref{kyr5}), we get $R_{i{\kern 1pt}
i} = {\rm{rdet}} _{{k}}\, {\rm {\bf A}}^{{i{\kern 1pt} i}}$, $k =
\min {\left\{ {I_{n} \setminus {\left\{ {i} \right\}}} \right\}}$.

Now suppose that $i \ne j$. Consider  monomials of
${\rm{rdet}}_{{i}}\, {\rm {\bf A}}$ in (\ref{kyr4})  such that
$a_{i{\kern 1pt} j} $ is the first from the left in each of their:
\[
\begin{array}{c}
   a_{i{\kern 1pt} j} \cdot R_{i{\kern 1pt} j} = {\sum\limits_{\bar
{\sigma} \in\,  S_{n}}  {\left( { - 1} \right)^{n - r}a_{i{\kern
1pt} j}\,  a_{j{\kern 1pt}  i_{k_{1}} }  \ldots } } a_{i_{k_{1} +
l_{1}}   i}  \ldots  a_{i_{k_{r}}  i_{k_{r} + 1}} \ldots
a_{i_{k_{r} + l_{r}}  i_{k_{r}} }  =\\
    = - a_{i{\kern 1pt} j} \cdot {\sum\limits_{\bar {\sigma}  \in \,S_{n}}
{\left( { - 1} \right)^{n - r - 1}a_{j {\kern 1pt} i_{k_{1}} }
\ldots } } a_{i_{k_{1} + l_{1}}   i} \ldots a_{i_{k_{r}}  i_{k_{r}
+ 1}}   \ldots a_{i_{k_{r} + l_{r}}  i_{k_{r}} } ,\\
\end{array}
\]
 \[ \bar {\sigma}  = \left( {i\,j\,\,i_{k_{1}}  \ldots i_{k_{1} +
l_{1}} } \right)\ldots \left( {i_{k_{r}}  i_{k_{r} + 1} \ldots
i_{k_{r} + l_{r}} } \right) \quad  {\forall r = \overline {1,n -
1}}.
\]
Denote $\tilde {a}_{i_{k_{1} + l_{1}}  j} = a_{i_{k_{1} + l_{1}}
i} $, $\left( {\forall}  i_{k_{1} + l_{1}}  \in  {I_{n}} \right)$.
Then
\[
\begin{array}{c}
   a_{i{\kern 1pt} j} \cdot R_{i{\kern 1pt} j} =
 - a_{i{\kern 1pt} j} \cdot {\sum\limits_{\bar {\sigma} _{1} \in \,
 \mathord{\buildrel{\lower3pt\hbox{$\scriptscriptstyle\frown$}}\over{S}_{n-1}}}
 {\left( { - 1} \right)^{n - r -
1}a_{j{\kern 1pt} i_{k_{1}} }   \ldots } } \tilde {a}_{i_{k_{1} +
l_{1} }  j}    \ldots  a_{i_{k_{r} + l_{r}}   i_{k_{r}} },\\
  \bar {\sigma} _{1} = \left( {j\,i_{k_{1}}  \ldots i_{k_{1} +
l_{1}} } \right)\ldots \left( {i_{k_{r}}  i_{k_{r} + 1} \ldots
i_{k_{r} + l_{r}} } \right) \quad  {\forall r = \overline {1,n -
1}}.
\end{array}
\]
The permutation   $\bar {\sigma} _{1} $  does not contain the
index $i$ in every monomial of $R_{i{\kern 1pt} j}$. By
(\ref{kyr4}) this permutation satisfies the conditions of
Definition \ref{kyrc1} for  ${\rm rdet}_{j} {\rm {\bf A}}_{.{\kern
1pt} j}^{i{\kern 1pt} i} \left( {{\rm {\bf a}}_{.\,{\kern 1pt}
{\kern 1pt} i}} \right)$. Here the matrix ${\rm {\bf A}}_{.{\kern
1pt} j}^{i{\kern 1pt} i} \left( {{\rm {\bf a}}_{.\,{\kern 1pt}
{\kern 1pt} i}} \right) $  is obtained from ${\rm {\bf A}}$ by
replacing the $j$th column with the column $i$, and then by
deleting both the $i$th row and  column. That is,
\[
{\sum\limits_{\bar {\sigma} _{1} \in \,
 \mathord{\buildrel{\lower3pt\hbox{$\scriptscriptstyle\frown$}}\over{S}_{n-1}}}
  {\left( { - 1} \right)^{n - r - 1}a_{j{\kern 1pt}
 i_{k_{1}} }  \ldots} }\tilde {a}_{i_{k_{1} + l_{1}}
{\kern 1pt} j} \ldots a_{i_{k_{r} + l_{r}}  {\kern 1pt} i_{k_{r}}
}  = {\rm{rdet}} _{{j}}\, {\rm {\bf A}}_{{.j}}^{{i{\kern 1pt}i}}
\left( {{\rm {\bf a}}_{{.{\kern 1pt} i}}} \right)
\]

Therefore, $R_{ij} = - {\rm rdet} _{j} {\rm {\bf A}}_{.{\kern 1pt}
j}^{i{\kern 1pt} i} \left( {{\rm {\bf a}}_{.{\kern 1pt} {\kern
1pt} i}} \right)$, if $i \ne j$.$\blacksquare$

\begin{definition} The $j$th column determinant
 of  ${\rm {\bf A}}\in {\rm M}\left( {n,\bf H}
\right)$
 is defined
as the alternating sum of  $n!$ products  of entries of ${\rm {\bf
A}}$, during which the index permutation of every product is
written by the direct product of disjoint cycles. If the
permutation is even, then
 products  of entries has a sign "$+$". If the permutation is
odd, then  products  of entries has a sign "$-$". That is
\[
{\rm{cdet}} _{{j}}\, {\rm {\bf A}} = {{\sum\limits_{\tau \in
S_{n}} {\left( { - 1} \right)^{n - r}a_{j_{k_{r}} j_{k_{r} +
l_{r}} } \ldots a_{j_{k_{r} + 1}  j_{k_{r}} }  \ldots } }a_{j\,
j_{k_{1} + l_{1}} }  \ldots  a_{ j_{k_{1} + 1} j_{k_{1}}
}a_{j_{k_{1}} j}},
\]
where $S_{n} $ is the symmetric group on the set
$J_{n}=\{1,\ldots,n\}$.  The  right-ordered cycle notation of the
permutation $\tau \in S_{n}$ is written as follows:
 \[
 \tau = \left( {j_{k_{r} + l_{r}}  \ldots j_{k_{r} + 1} j_{k_{r}}
} \right)\ldots \left( {j_{k_{2} + l_{2}}  \ldots j_{k_{2} + 1}
j_{k_{2}} } \right){\kern 1pt} \left( {j_{k_{1} + l_{1}}  \ldots
j_{k_{1} + 1} j_{k_{1} } j} \right).
\]
Here the index $j$ opens the first cycle from the right  and other
cycles satisfy the following conditions
\[
 j_{k_{2}}  < j_{k_{3}}  < \ldots < j_{k_{r}},
\quad j_{k_{t}}  < j_{k_{t} + s}, \quad \left( {\forall t =
\overline {2,r}}  \right), \quad \left( {\forall s = \overline
{1,l_{t}} }  \right).
\]
\end{definition}
\newtheorem{remark}{Remark}[section]
\begin{remark}
A peculiarity of calculation of column determinants is such that
coefficients of every monomials are written from right to left.
\end{remark}
\begin{lemma} Let $L_{i{\kern 1pt} j} $ be
the left $ij$-th cofactor of  a matrix ${\rm {\bf A}}\in {\rm
M}\left( {n,\bf H} \right)$, that is $ {\rm{cdet}} _{{j}}\, {\rm
{\bf A}} = {{\sum\limits_{i = 1}^{n} {L_{i{\kern 1pt} j} \cdot
a_{i{\kern 1pt} j}} }}$, $(\forall {j = \overline {1,n}} )$. Then
\[
L_{i{\kern 1pt} j} = {\left\{ {\begin{array}{*{20}c}
 -{\rm{cdet}} _{i}\, {\rm {\bf A}}_{i{\kern 1pt} .}^{j{\kern 1pt}j}
  \left( {{\rm {\bf a}}_{j{\kern 1pt}. } }\right),& {i \ne
j},\\
 {\rm{cdet}} _{k}\, {\rm {\bf A}}^{j\, j},& {i = j},
\\
\end{array} }\right.}
\]
where  ${\rm {\bf A}}_{i{\kern 1pt} .}^{jj} \left( {{\rm {\bf
a}}_{j{\kern 1pt} .} } \right)$ is obtained from ${\rm {\bf A}}$
 by replacing the $i$th row with the $j$th row, and then by
deleting both the $j$th row and  column; $k = \min {\left\{
{J_{n}} \right.} \setminus {\left. {\{j\}} \right\}} $.
\end{lemma}
The proof of this lemma is similar to that of Lemma \ref{kyrc2}.
\begin{remark} Clearly, any monomial of  each row
or column determinant of a square matrix corresponds to a certain
monomial of  another row or column determinant such that both of
them consist of the same coefficients and vary  in their ordering
only. If  the entries of ${\rm {\bf A}}$ are commutative, then $
{\rm{rdet}} _{1}\, {\rm {\bf A}} = \ldots = {\rm{rdet}} _{{n}}
{\rm {\bf A}} = {\rm{cdet}} _{1}\, {\rm {\bf A}} = \ldots =
{\rm{cdet}} _{{n}} {\rm {\bf A}}. $
\end{remark}

 Consider
the basic properties of the column and row determinants of a
square matrix over $\bf H$,  the proofs of which   immediately
follow from the definitions.
\newtheorem{theorem}{Theorem}[section]
\begin{theorem}
If one of the rows (columns) of  ${\rm {\bf A}}\in {\rm M}\left(
{n,\bf H} \right)$  consists of zeros only, then $
{\rm{rdet}}_{{i}}\, {\rm {\bf A}} = 0$ and ${\rm{cdet}} _{{i}}\,
{\rm {\bf A}} = 0$, $ (\forall {i = \overline {1,n}} ). $
\end{theorem}
\begin{theorem} If the $i$th row of  ${\rm {\bf A}}\in {\rm M}\left( {n,\bf H}
\right)$ is left-multiplied by $b \in \bf H $, then $
{\rm{rdet}}_{{i}}\, {\rm {\bf A}}_{{i{\kern 1pt}.}} \left( {b
\cdot {\rm {\bf a}}_{{i{\kern 1pt}.}}} \right) = b \cdot
{\rm{rdet}}_{{i}}\, {\rm {\bf A}}$, $ (\forall {i = \overline
{1,n}} ).$
\end{theorem}
\begin{theorem} If the $j$th column of
 ${\rm {\bf A}}\in {\rm M}\left( {n,\bf H}
\right)$ is right-multiplied by $b \in \bf H$, then $ {\rm{cdet}}
_{{j}}\, {\rm {\bf A}}_{{.{\kern 1pt}j}} \left( {{\rm {\bf
a}}_{{.{\kern 1pt}j}} \cdot b} \right) = {\rm{cdet}} _{{j}}\, {\rm
{\bf A}}$, $(\forall {j = \overline {1,n}} ).$
\end{theorem}
\begin{theorem}\label{kyrc4} If for  ${\rm {\bf A}}\in {\rm M}\left( {n,\bf H}
\right)$\, $\exists t \in I_{n} $ such that $a_{tj} = b_{j} +
c_{j} $\, $\left( {\forall j = \overline {1,n}} \right)$, then
\[
\begin{array}{l}
   {\rm{rdet}}_{{i}}\, {\rm {\bf A}} = {\rm{rdet}}_{{i}}\, {\rm {\bf
A}}_{{t{\kern 1pt}.}} \left( {{\rm {\bf b}}} \right) +
{\rm{rdet}}_{{i}}\, {\rm {\bf A}}_{{t{\kern 1pt}.}} \left( {{\rm
{\bf c}}} \right), \\
  {\rm{cdet}} _{{i}}\, {\rm {\bf A}} = {\rm{cdet}} _{{i}}\, {\rm
{\bf A}}_{{t{\kern 1pt}.}} \left( {{\rm {\bf b}}} \right) +
{\rm{cdet}}_{{i}}\, {\rm {\bf A}}_{{t{\kern 1pt}.}} \left( {{\rm
{\bf c}}} \right) \quad  {\forall {i = \overline {1,n}}},
\end{array}
\]
where ${\rm {\bf b}}=(b_{1},\ldots, b_{n})$, ${\rm {\bf
c}}=(c_{1},\ldots, c_{n}).$
\end{theorem}
\begin{theorem} \label{kyrc5} If for ${\rm {\bf A}}\in {\rm M}\left( {n,\bf H}
\right)$\, $\exists t \in J_{n} $ such that $a_{i\,t} = b_{i} +
c_{i}$ $\left( {\forall i = \overline {1,n}} \right)$, then
\[
\begin{array}{l}
  {\rm{rdet}}_{{j}}\, {\rm {\bf A}} = {\rm{rdet}}_{{j}}\, {\rm {\bf
A}}_{{.\,{\kern 1pt}t}} \left( {{\rm {\bf b}}} \right) +
{\rm{rdet}}_{{j}}\, {\rm {\bf A}}_{{.\,{\kern 1pt} t}} \left(
{{\rm
{\bf c}}} \right),\\
  {\rm{cdet}} _{{j}}\, {\rm {\bf A}} = {\rm{cdet}} _{{j}}\, {\rm
{\bf A}}_{{.\,{\kern 1pt}t}} \left( {{\rm {\bf b}}} \right) +
{\rm{cdet}} _{{j}} {\rm {\bf A}}_{{.\,{\kern 1pt}t}} \left( {{\rm
{\bf c}}} \right) \quad  {\forall {j = \overline {1,n}}},
\end{array}
\]
where ${\rm {\bf b}}=(b_{1},\ldots, b_{n})^T$, ${\rm {\bf
c}}=(c_{1},\ldots, c_{n})^T.$
\end{theorem}
\begin{theorem}  If ${\rm {\bf A}}^{ *} $
is the Hermitian adjoint matrix of  ${\rm {\bf A}}\in {\rm
M}\left( {n,\bf H} \right)$, then $ {\rm{rdet}}_{{i}}\, {\rm {\bf
A}}^{ *} =  \overline{{{\rm{cdet}} _{{i}}\, {\rm {\bf A}}}}$,
$\left( {\forall {i = \overline {1,n}}} \right). $
\end{theorem}
\begin{remark} Since the column and row determinants of an arbitrary square
matrix over the quaternion skew field do not satisfy Axiom 1 but
these determinants are defined by analogy to the determinant of
 a complex square matrix, then we consider theirs  as pre-determinants.
\end{remark}
\section{ A determinant of a Hermitian matrix}
The following lemma is needed for the sequel.
\begin{lemma} \label{kyrc3}  Let $T_{n} $ be the sum
of all  possible products  of  the $n$ factors, each of which are
either $h_{i} \in {\bf H}$ or $ \overline {h_{i}} $, $(\forall i =
\overline {1,n}), $ by specifying the ordering in the terms, i.e.:
\[
T_{n} = h_{1} \cdot h_{2} \cdot \ldots \cdot h_{n} + \overline
{h_{1}} \cdot h_{2} \cdot \ldots \cdot h_{n} + \ldots + \overline
{h_{1}}  \cdot \overline {h_{2}}  \cdot \ldots \cdot \overline
{h_{n}}  {\rm .}
\]
Then $T_{n}$ consists of the  $2^{n}$ terms and $T_{n} = {\rm
t}\left( {h_{1}} \right)\;{\rm t}\left( {h_{2}} \right)\;\ldots
\;{\rm t}\left( {h_{n}} \right).$
\end{lemma}
{\textit{Proof.}} The number $2^{n}$ of terms of the sum $T_{n}$
is equal to the number of  ordered  combinations of $n$ unknown
elements with two values.

  The proof goes by induction on $n$.

\noindent (i) If $n = 1$, then $T_{1} = \overline {h_{1}}  + h_{1}
= {\rm t}\left( {h_{1}} \right)$.

\noindent (ii) Suppose the lemma is true for $n - 1$:
\[
\begin{array}{c}
 T_{n - 1} = h_{1} \cdot h_{2} \cdot \ldots \cdot h_{n - 1} +
\overline {h_{1}}  \cdot h_{2} \cdot \ldots \cdot h_{n - 1} +
\ldots + \overline {h_{1}}  \cdot \overline {h_{2}}  \cdot \ldots
\cdot \overline {h_{n - 1}} =\\
 = {\rm t}\left( {h_{1}}  \right)\;{\rm t}\left( {h_{2}}  \right)\;\ldots \;{\rm t}\left(
{h_{n - 1}}  \right).
\end{array}
\]
\noindent (iii) Now we prove
 that it is valid  for $n$.
\[
T_{n} = h_{1} \cdot h_{2} \cdot \ldots \cdot h_{n} + \overline
{h_{1}} \cdot h_{2} \cdot \ldots \cdot h_{n} + \ldots + \overline
{h_{1}}  \cdot \overline {h_{2}}  \cdot \ldots \cdot \overline
{h_{n}}.
\]
By factoring  the  common right-side factors either $\overline
{h_{n}} $ or $h_{n} $ respectively, we obtain
\[
\begin{array}{c}
  T_{n} = T_{n - 1} \cdot h_{n} + T_{n - 1} \cdot \overline {h_{n}}
= T_{n - 1} \cdot \left( {h_{n} + \overline {h_{n}} }  \right) =
T_{n - 1} \cdot {\rm t}\left( {h_{n}}  \right) = \\
   = {\rm t}\left( {h_{1}}  \right) \cdot \;{\rm t}\left( {h_{2}}  \right) \cdot \;\ldots
\cdot {\rm t}\left( {h_{n - 1}}  \right) \cdot \;{\rm t}\left(
{h_{n}} \right).\blacksquare
\end{array}
\]
\begin{theorem} If ${\rm {\bf A}}\in {\rm M}\left(
{n,\bf H} \right)$ is a Hermitian matrix, then
\[
{\rm{rdet}} _{1} {\rm {\bf A}} = \ldots = {\rm{rdet}} _{n} {\rm
{\bf A}} = {\rm{cdet}} _{1} {\rm {\bf A}} = \ldots = {\rm{cdet}}
_{n} {\rm {\bf A}}  \in {\bf R}.
\]
\end{theorem}
{\textit{Proof.}} At first we note  that if a matrix ${\rm {\bf
A}}=(a_{ij})\in {\rm M}\left( {n,\bf H} \right)$ is  Hermitian,
then we have $a_{ii} \in{\bf R}$ and $a_{ij} = \overline
{a_{ji}}$, $(\forall i,j = \overline {1,n} )$.
 We divide the set of monomials of some
${\rm{rdet}}_{i} {\rm {\bf A}}
 \left( {\forall i = \overline {1,n}}  \right)$ into two
subsets. If indices of coefficients of monomials form permutations
as products of disjoint cycles of length 1 and 2, then we include
these monomials in the first subset. Other monomials belong to the
second subset. If indices of coefficients form a disjoint cycle of
length 1, then these coefficients are entries of the principal
diagonal of the Hermitian matrix ${\rm {\bf A}}$. Hence, they are
real numbers. If indices of coefficients form a disjoint cycle of
length 2, then these elements are conjugated, $a_{i_{k} i_{k + 1}}
= \overline {a_{i_{k + 1} i_{k}} } $, and their product  takes on
a real value as well,
\[
a_{i_{k} i_{k + 1}}  \cdot a_{i_{k + 1} i_{k}}  = \overline
{a_{i_{k + 1} i_{k}} }  \cdot a_{i_{k + 1} i_{k}}  = {\rm
n}(a_{i_{k + 1} i_{k}}  ) \in {\bf R}.
\]
So, all monomials of the first subset take on real values.

Now we consider some monomial  $d$ from the second subset. Assume
that  indices of its coefficients form a permutation as a product
of $r$ disjoint cycles. Denote $i_{k_{1}}:=i$.
\begin{equation}
\label{kyr6}
\begin{array}{l}
 d = ( - 1)^{n - r}a_{i_{k_{1}}i_{k_{1}+1} }   \ldots  a_{i_{k_{1} + l_{1}
} i_{k_{1}}}  a_{i_{k_{2}}  i_{k_{2} + 1}}   \ldots  a_{i_{k_{2} +
l_{2}} i_{k_{2}} }   \ldots  a_{i_{k_{m}}  i_{k_{m} + 1}} \ldots
\times
\\
 \times a_{i_{k_{m} + l_{m}}  i_{k_{m}} }   \ldots  a_{i_{k_{r}}
i_{k_{r} + 1}}  \ldots  a_{i_{k_{r} + l_{r}}  i_{k_{r}} }  = ( -
1)^{n - r}h_{1}  h_{2}  \ldots  h_{m}  \ldots  h_{r} ,
\\
 \end{array}
\end{equation}
\noindent where $h_{s} = a_{i_{k_{s}} i_{k_{s} + 1}}  \cdot \ldots
\cdot a_{i_{k_{s} + l_{s}}  i_{k_{s}} }$, $\left( {\forall s =
\overline {1,r}} \right)$, and $ m \in \{1,\ldots ,r\}.$   If
$l_{s} = 1$, then $h_{s} = a_{i_{k_{s}}  i_{k_{s} + 1}}
a_{i_{k_{s} + 1{\kern 1pt}}  i_{k_{s}} }  = {\rm n}(a_{i_{k_{s}}
i_{k_{s} + 1}} ) \in {\bf R}$. If $l_{s} = 0$, then $h_{s} =
a_{i_{k_{s}} i_{k_{s}} }  \in {\bf R}$. If $l_{s} = \overline
{0,1}$ for $\forall s = \overline {1,r} $  in (\ref{kyr6}), then
we obtain a monomial of the first subset. Let $\exists s \in I_{n}
$ such that $l_{s} \ge 2$. Then
\[
   \overline {h_{s}}  = \overline {a_{i_{k_{s}}  i_{k_{s} + 1}}
\ldots  a_{i_{k_{s} + l_{s}}  i_{k_{s}} } }
   = \overline
{a_{i_{k_{s} + l_{s} } i_{k_{s}} } }   \ldots  \overline
{a_{i_{k_{s}}  i_{k_{s} + 1}} } = a_{i_{k_{s}}  i_{k_{s} + l_{s}}
}  \ldots  a_{i_{k_{s} + 1} i_{k_{s}} }.
\]
Denote by  $\sigma _{s} \left( {i_{k_{s}} }  \right){\rm :} =
\left( {i_{k_{s}}  i_{k_{s} + 1} \ldots i_{k_{s} + l_{s}} }
\right)$
  a disjoint cycle of indices of $d$, $(\forall s =
\overline {1,r})$. The disjoint cycle $\sigma _{s} \left(
{i_{k_{s}} }  \right)$ corresponds to the factor $h_{s}$. Then
$\sigma _{s}^{ - 1} \left( {i_{k_{s}} }  \right) = \left(
{i_{k_{s}} i_{k_{s} + l_{s}} i_{k_{s} + 1} \ldots i_{k_{s} + 1}}
\right)$ is the inverse disjoint cycle and $\sigma _{s}^{ - 1}
\left( {i_{k_{s}} } \right) $ corresponds to $\overline{h_{s}}$,
$(\forall s = \overline {1,r})$. By Lemma \ref{kyrc3} for $d$
there exist another $2^{p}-1$ monomials, (where $p = r - \rho $
and $\rho $ is the number of disjoint cycles of length 1 and 2),
such that their index permutations are written by the direct
products of $r$ disjoint cycles either $\sigma _{s} \left(
{i_{k_{s}} } \right)$ or $\sigma _{s}^{ - 1} \left( {i_{k_{s}} }
\right)$ by specifying their ordering by $s$ from $1$ to $r$.
These permutations are left ordered in the cycle representation
according to Definition \ref{kyrc1}. Suppose $C_{1} $ is the sum
of these $2^{p}-1$ monomials and $d$. Then by Lemma \ref{kyrc3} we
obtain
\[
C_{1} = ( - 1)^{n - r}\alpha \;{\rm t}(h_{\nu _{1}}  )\;\ldots
\;{\rm t}(h_{\nu _{p}}  )  \in {\bf R}.
\]
\noindent Here $\alpha \in {\bf R}$ is the product of coefficients
 whose indices   form disjoint cycles of length 1 and 2,
  $\nu _{k} \in \{1,\ldots ,r\}$, $(\forall k = \overline {1,p} ).$
Thus for an arbitrary monomial from the second subset of
${\rm{rdet}}_{i}\, {\rm {\bf A}}$, we can  find the $2^{p}$
monomials  such that their sum  takes on a real value. Therefore,
${\rm{rdet}} _{i}\, {\rm {\bf A}}\in {\bf R}$.

Now we  prove the equality of all  row determinants of ${\rm {\bf
A}}$. Consider  ${\rm{rdet}} _{j}\, {\rm {\bf A}}$ such that $j
\ne i, \, \left( {\forall j = \overline {1,n}} \right)$. We divide
the set of monomials of ${\rm{rdet}} _{j}\, {\rm {\bf A}}$ into
two subsets using the same rule as for ${\rm{rdet}}_{i}\, {\rm
{\bf A}}.$ Monomials from the first subset are products of the
real factors (either entries of the principal diagonal of ${\rm
{\bf A}}$ or  norms of entries). Hence each monomial from the
first subset of ${\rm{rdet}}_{i}\, {\rm {\bf A}}$ is equal to a
corresponding monomial from the first subset of ${\rm{rdet}}
_{j}\, {\rm {\bf A}}$. Now consider the monomial $d_{1} $ from the
second subset of monomials of ${\rm{rdet}}_{i}\, {\rm {\bf A}}$
consisting of coefficients that are equal to the coefficients of
$d$ but are
 placed in another arrangement.
 Consider all  possibilities
of the arrangement of coefficients in $d_{1} $.

 (i) Suppose indices of its coefficients form a
permutation as a product of $r$ disjoint cycles and these cycles
coincide with the $r$ disjoint cycles of $d$. But the index
permutation  of $d$ distinguishes from the index  permutation of
$d_{1}$ by the ordering of disjoint cycles. Then we have
\[d_{1} = ( - 1)^{n - r}\alpha h_{\mu}   \ldots  h_{\lambda}  ,
\]
where $\{\mu ,\ldots ,\lambda \} = \{\nu _{1} ,\ldots ,\nu _{p}
\}$. By Lemma \ref{kyrc3} there exist $2^{p }- 1$ monomials among
the monomials from the second subset of ${\rm{rdet}} _{j}\, {\rm
{\bf A}}$ such that each of them is equal to a product of $p$
factors either $h_{s} $ or $\overline {h_{s}} $, $(\forall s \in
\{\mu ,\ldots ,\lambda \})$, multiplied by $( - 1)^{n - r}\alpha
$. Hence by Lemma \ref{kyrc3}, we obtain
\[
C_{2} = ( - 1)^{n - r}\alpha \;t(h_{\mu}  )\;\ldots
\;t(h_{\lambda}  ) = ( - 1)^{n - r}\;\alpha \;t(h_{\nu _{1}}
)\ldots \;t(h_{\nu _{p}}  ) = C_{1} .
\]

 (ii) Now suppose that  in addition to the case (i)  the index $j$
 is
placed inside some disjoint cycle of  indices of $d$, e.g. $ j \in
\{i_{k_{m+1}},...,i_{k_{m}+l_{m}}\}$. Denote $j = i_{k_{m} + q} $.
Then $d_{1} $ is represented as follows:
\begin{equation}\label{kyr7}
\begin{array}{c}
   d_{1} = ( - 1)^{n - r}a_{i_{k_{m} + q} i_{k_{m} + q + 1}}  \ldots
\quad  a_{i_{k_{m} + l_{m}}  i_{k_{m}} }  \, a_{i_{k_{m}} i_{k_{m}
+ 1}} \ldots   \times \\
   \times   a_{i_{k_{m} + q - 1} i_{k_{m} + q}} a_{i_{k_{\mu} }  i_{k_{\mu}  + 1}}
      \ldots  a_{i_{k_{\mu
} + l_{\mu} }  i_{k_{\mu} } }   \ldots a_{i_{k_{\lambda} }
i_{k_{\lambda}  + 1}}  \ldots a_{i_{k_{\lambda}  + l_{\lambda} }
i_{k_{\lambda} } }  = \\
  =(-1)^{n - r}\alpha  \tilde {h}_{m}  h_{\mu}
\ldots  h_{\lambda},
\end{array}
\end{equation}
\noindent where $\{m,\mu ,\ldots ,\lambda \} = \{\nu _{1} ,\ldots
,\nu _{p} \}$. Except for $\tilde {h}_{m} $,  each factor of
$d_{1}$ in (\ref{kyr7}) corresponds to the equal factor of $d$ in
(\ref{kyr6}).  We have $t(\tilde {h}_{m} ) = t(h_{m} )$ by the
rearrangement property of the trace. Hence by Lemma \ref{kyrc3}
and by analogy to the previous case, we obtain the following
equality.
\[
\begin{array}{c}
  C_{2} = ( - 1)^{n - r}\alpha \;t(\tilde {h}_{m} )\;t(h_{\mu}
)\;\ldots \;t(h_{\lambda}  )=\\
 = ( - 1)^{n - r}\;\alpha \;t(h_{\nu _{1}}  )\;\ldots \;t(h_{m} )\;\ldots
\;t(h_{\nu _{p}}  ) = C_{1} .
\end{array}
\]

(iii) If in addition to the case (i)  the index $i$ is placed
inside some disjoint cycles of the index permutation of $d_{1}$,
then we apply the rearrangement property of the trace to a factor
whose indices belong to this cycle. As in the previous cases among
monomials from the second subset of ${\rm{rdet}} _{j} \,{\rm {\bf
A}}$, we find  $2^{p}$ monomials such that by Lemma \ref{kyrc3}
their sum is equal to the sum of the corresponding $2^{p}$
monomials of ${\rm{rdet}}_{i} {\rm {\bf A}}$. Clearly, we obtain
the same conclusion at association of all previous cases, then we
apply  the rearrangement property of the trace twice.

Thus, in any case  each  sum of  $2^{p}$ corresponding monomials
from the second subset of ${\rm{rdet}} _{{j}}\, {\rm {\bf A}}$ is
equal to the sum of  $2^{p}$ monomials of ${\rm{rdet}}_{{i}}\,
{\rm {\bf A}}$. Here $p$ is the number of disjoint cycles of
length  more than 2. Thus,
\[
{\rm{rdet}}_{{i}}\, {\rm {\bf A}} = {\rm{rdet}} _{{j}}\, {\rm {\bf
A}}  \in {\bf R},\quad \forall i, j = \overline {1,n}.
\]

Now we prove the equality ${\rm{cdet}} _{{i}}\, {\rm {\bf A}} =
{\rm{rdet}}_{{i}}\, {\rm {\bf A}}$\, $\left( {\forall i =
\overline {1,n}}  \right)$. We divide the set of monomials of
${\rm{cdet}} _{{i}}\, {\rm {\bf A}}$ into two subsets conforming
to  the same rule as for ${\rm{rdet}}_{{i}}\, {\rm {\bf A}}$. Each
monomial from the first subset of ${\rm{cdet}} _{{i}}\, {\rm {\bf
A}}$ is equal to the corresponding monomial of
${\rm{rdet}}_{{i}}\, {\rm {\bf A}}$ since their factors are real
numbers (either entries of the principal diagonal of ${\rm {\bf
A}}$ or  norms of entries of ${\rm {\bf A}}$). Consider the
monomial $d_{2} $ from the second subset of monomials of
${\rm{cdet}}_{i}\, {\rm {\bf A}}$  consisting of coefficients that
are equal to the coefficients of $d$. The coefficients of $d_{2} $
are
 placed in the same ordering as for $d$ but from left to right.
 If $\rho $ is the number of disjoint cycles of length 1 and 2,
 and
  $p = r - \rho $, then
\[
\begin{array}{c}
 d_{2} = ( - 1)^{n - r}a_{i_{k_{r}}  i_{k_{r} + l_{r}} }   \ldots
a_{i_{k_{r} + 1} i_{k_{r}} }   \ldots  a_{i_{k_{2}}  i_{k_{2} +
l_{2}} }   \ldots  a_{i_{k_{2} + 1} i_{k_{2}} }\times\\
 \times a_{i_{k_{1}}i_{k_{1} + l_{1}} }   \ldots  a_{i_{k_{1}+1}  i_{k_{1}}} = ( -
1)^{n - r}\alpha \;h_{\tau _{p}}\ldots  h_{\tau _{1}}.
\end{array}
\]
Here $\alpha $ is a product of coefficients whose indices form
disjoint cycles of length 1 and 2. We have
\[
h_{\tau _{s}}  = a_{i_{k_{s}}  i_{k_{s} + l_{s}} }  \cdot \ldots
\cdot a_{i_{k_{s} + 1} i_{k_{s}} }  = \overline {a_{i_{k_{s}}
i_{k_{s} + 1}} \cdot \ldots \cdot a_{i_{k_{s} + l_{s}}  i_{k_{s}}
} }= \overline {h_{\nu _{s}} } \quad \forall s = \overline {1,p}.
\]
By Lemma \ref{kyrc3} among  monomials from the second subset of
${\rm{cdet}} _{{i}}\, {\rm {\bf A}}$,  there exist  $2^{p} - 1$
monomials for $d_{2}$ such  that each of them  is equal to a
product of $p$ factors, either $h_{\tau _{s}}$ or $\overline
{h_{\tau _{s}}} $\, $(s = \overline {1,p} )$, by specifying their
right-ordering, and is multiplied by $( - 1)^{n - r}\alpha $.
Consider the sum $C_{3}$ of these monomials and $d$. Due to
commutativity of real numbers and by Lemma \ref{kyrc3}, we get
\[
\begin{array}{c}
 C_{3} = ( - 1)^{n - r}\alpha \;t(h_{\tau _{p}}  )\;\ldots \;t(h_{\tau _{1}
} ) = ( - 1)^{n - r}\alpha \;t(\overline {h_{\nu _{p}} } )\;\ldots
\;t(\overline {h_{\nu _{1}} }  ) = \\
  = ( - 1)^{n - r}\alpha \;t(h_{\nu _{1}}  )\ldots \;t(h_{\nu _{p}}
) = C_{1} \\
 \end{array}
\]

Therefore,  each sum of the $2^{p}$ corresponding monomials from
the second subset of ${\rm{cdet}} _{{i}}\, {\rm {\bf A}}$ is equal
to a sum of the $2^{p}$ monomials of ${\rm{rdet}}_{{i}}\, {\rm
{\bf A}}$ and vice versa.

Thus, ${\rm{cdet}} _{{i}}\, {\rm {\bf A}} = {\rm{rdet}}_{{i}}\,
{\rm {\bf A}}  \in \bf R$ \, $\left( {\forall i = \overline {1,n}}
\right)$.$\blacksquare$

\begin{remark}
 Since all  column and  row determinants of a
Hermitian matrix over $\bf H$ are equal, we can define the
 determinant of a  Hermitian matrix ${\rm {\bf A}}\in {\rm M}\left(
{n,\bf H} \right)$. By definition, put
\[\det {\rm {\bf A}}: = {\rm{rdet}}_{{i}}\,
{\rm {\bf A}} = {\rm{cdet}} _{{i}}\, {\rm {\bf A}}, \quad {\forall
i = \overline {1,n}}.\] The determinant of a Hermitian matrix
satisfies Axiom 1. It follows from Theorem \ref{kyrc9} and
Corollary \ref{kyrc7}.

\end{remark}

\begin{remark} By Lemma \ref{kyrc3} we have
\begin{equation}\label{kyr8}
\det {\rm {\bf A}} = - {\sum\limits_{\sigma \in I_{n}} {a_{i{\kern
1pt} j} \cdot {\rm{rdet}} _{{j}}\, {\rm {\bf A}}_{.j}^{i{\kern
1pt}i} \left( {{\rm {\bf a}}_{.i}} \right) + a_{i{\kern 1pt} i}
\cdot {\rm{rdet}} _{{k}}\, {\rm {\bf A}}_{}^{i{\kern 1pt}i}} },\,k
= \min {\left\{ {I_{n}}  \right.} \setminus {\left. {\{i\}}
\right\}}.
\end{equation}
By comparing expressions (\ref{kyr1}) and (\ref{kyr8}) for a
Hermitian matrix ${\rm {\bf A}}\in {\rm M}\left( {n,\bf H}
\right)$, we conclude that the row determinant of a Hermitian
matrix coincides with the Moore determinant. Hence the row and
column determinants extend the Moore determinant to an arbitrary
square matrix.
\end{remark}

\section{ Properties of the column and  row \\ determinants of a Hermitian matrix}
\begin{theorem} \label{kyrc6} If the matrix ${\rm {\bf A}}_{j.} \left( {{\rm {\bf a}}_{i.}}
\right)$ is obtained from a Hermitian matrix ${\rm {\bf A}}\in
M\left( {n,\bf H} \right)$  by replacing  its $j$th row with the
$i$th row,  then
 \[{\rm{rdet}} _{{j}} {\rm {\bf A}}_{{j\,.}}
\left( {{\rm {\bf a}}_{{i\,.}}} \right) = 0,\,(\forall
i,j=\overline{1,n},\,i\neq j).\]
\end{theorem}
{\textit{Proof.}} We assume $n > 3$ for  ${\rm {\bf A}}\in {\rm
M}\left( {n,\bf H}\right)$. The case $n \le 3$ can be easily
proved by direct calculation. Consider some monomial $d$ of
${\rm{rdet}} _{{j}}\, {\rm {\bf A}}_{{j\,.}} \left( {{\rm {\bf
a}}_{{i\,.}}} \right)$. Suppose indices of its coefficients form a
permutation as a product of $r$ disjoint cycles, and  denote $i =
i_{s} $. Consider all
 possibilities of disposition of an entry of the $i_{s} $th
row  in the monomial $d$.

(i) Suppose an entry of the $i_{s} $th row is placed in $d$ such
that the index $i_{s} $ opens some disjoint cycle, i.e.:
\begin{equation}
\label{kyr9} d = ( - 1)^{n - r}a_{j{\kern 1pt} i_{1}}   \ldots
 a_{i_{k} j} \, u_{1} \ldots  u_{\rho}  \, a_{i_{s} i_{s +
1}}   \ldots  a_{i_{s + m} i_{s}}  \, v_{1} \ldots v_{p}
\end{equation}
\noindent Here by $u_{\tau}  $ and $v_{t} $ we denote   products
of coefficients whose indices form some disjoint cycles $(\forall
\tau = \overline {1,\rho} ,\;\forall t = \overline {1,p},\, \rho +
p = r - 2)$ or there are no such products. For $d$ there are the
following three monomials of ${\rm{rdet}} _{{j}}\, {\rm {\bf
A}}_{{j\,.}} \left( {{\rm {\bf a}}_{i\,.}} \right)$.
\[\begin{array}{c}
   d_{1} = ( - 1)^{n - r + 1}a_{j{\kern 1pt} i_{s + 1}} \ldots
 a_{i_{s + m} i_{s}}  \, a_{i_{s} i_{1}}   \ldots
a_{i_{k} j} \, u_{1}  \ldots u_{\rho}\, v_{1}
 \ldots v_{p},\\
   d_{2} = ( - 1)^{n - r + 1}a_{j{\kern 1pt} i_{s + m}}   \ldots
 a_{i_{s + 1} i_{s}}  a_{i_{s}{\kern 1pt} i_{1}} \ldots
 a_{i_{k} j} \, u_{1} \ldots  u_{\rho}  \, v_{1} \ldots
v_{p},\\
  d_{3} = ( - 1)^{n - r}a_{j{\kern 1pt} i_{1}}   \ldots a_{i_{k} j}
\, u_{1}  \ldots  u_{\rho}  \, a_{i_{s} i_{s + m}} \ldots  a_{i_{s
+ 1} i_{s}}\, v_{1} \ldots v_{p}.
\end{array}
\]

If $a_{j i_{1}}  \ldots  a_{i_{k} j} = x$ and $a_{i_{s} i_{s + 1}}
\ldots a_{i_{s + m} i_{s}}  = y$, then $\overline {y} = a_{i_{s}
i_{s + m}}  \ldots a_{i_{s + 1} \,i_{s}} $. Taking into account
$a_{j{\kern 1pt} i_{1}}  = a_{i_{s} i_{1}}  $, $a_{j{\kern
1pt}i_{s - 1}}  = a_{i_{s} i_{s - 1}}$ and $ a_{j{\kern 1pt}i_{s +
1}} = a_{i_{s} i_{s + 1}}  $, we consider the sum of these
monomials.
\[
\begin{array}{c}
d + d_{1} + d_{2} + d_{3} = ( - 1)^{n - r}(x  u_{1} \ldots
u_{\rho} \, y - y  x  u_{1} \ldots u_{\rho}- \overline {y} \cdot x
u_{1} \ldots  u_{\rho} +\\
  + x  u_{1}  \ldots u_{\rho} \overline {y}
)  v_{1}  \ldots v_{p} =( - 1)^{n - r} ( xu_{1} \ldots u_{\rho}
t(y)-t(y)xu_{1} \ldots u_{\rho}) v_{1}  \ldots  v_{p} =0
\end{array}
\]
Thus among the monomials of ${\rm{rdet}} _{{j}}\, {\rm {\bf
A}}_{{j\,.}} \left( {{\rm {\bf a}}_{{i\,.}}}  \right)$  we find
 three monomials  for $d$ such  that the sum of these  monomials and
 $d$ is equal to zero.

If in (9) $m = 0$ or $m = 1$, we  obtain such monomials
accordingly:
\[\begin{array}{c}
   \tilde {d} = ( - 1)^{n - r}a_{j{\kern 1pt} i_{1}}  \ldots
\cdot a_{i_{k} j} \, u_{1} \ldots  u_{\rho}  \, a_{i_{s} i_{s}} \,
v_{1} \ldots  v_{p},\\
    \mathord{\buildrel{\lower3pt\hbox{$\scriptscriptstyle\frown$}}\over {d}} = (
- 1)^{n - r}a_{j{\kern 1pt} i_{1}}   \ldots  a_{i_{k} j} \, u_{1}
\ldots  u_{\rho}  \, a_{i_{s} i_{s + 1}} \, a_{i_{s + 1} i_{s}} \,
v_{1}  \ldots  v_{p}.
\end{array}
   \]
There are the following  monomials for them:
\[\begin{array}{c}
     \tilde {d}_{1} = ( - 1)^{n - r + 1}a_{j{\kern 1pt}i_{s}}  \, a_{i_{s}
{\kern 1pt} i_{1}}  \, \ldots  a_{i_{k} j} \, u_{1}
 \ldots  u_{\rho}  \, v_{1}  \ldots v_{p},\\
   \mathord{\buildrel{\lower3pt\hbox{$\scriptscriptstyle\frown$}}\over {d}}
_{1} = ( - 1)^{n - r + 1}a_{j{\kern 1pt}i_{s + 1}} \, a_{i_{s + 1}
{\kern 1pt} i_{s}}  \,a_{i_{s} {\kern 1pt} i_{1}} \ldots a_{i_{k}
j}\, u_{1} \ldots  u_{\rho}\,
 v_{1}  \ldots  v_{p}.
\end{array}
 \]
Taking into account $a_{j{\kern 1pt} i_{1}}  = a_{i_{s} i_{1}},$
$a_{j{\kern 1pt} i_{s}}  = a_{i_{s} i_{s}}\in {\bf R}$,
$a_{j{\kern 1pt} i_{s + 1}} = a_{i_{s} i_{s + 1}}  $, and
 $ {a_{i_{s} i_{s + 1}} a_{i_{s + 1} i_{s}}}
\in {\bf R}$, we get $\tilde {d} + \tilde {d}_{1} = 0,$ \,
$\mathord{\buildrel{\lower3pt\hbox{$\scriptscriptstyle\frown$}}\over
{d}} +
\mathord{\buildrel{\lower3pt\hbox{$\scriptscriptstyle\frown$}}\over
{d}} _{1} = 0$. Hence,  the sums of corresponding two monomials of
${\rm{rdet}} _{{j}}\, {\rm {\bf A}}_{{j\,.}} ( {{\rm {\bf
a}}_{{i\,.}}})$ are equal to zero in this case.

 ii) Now suppose that the index $i_{s} $  is placed in  another disjoint
  cycle  than the index $j$ and does not open this cycle,
\[
  \mathord{\buildrel{\lower3pt\hbox{$\scriptscriptstyle\smile$}}\over
{d}} = ( - 1)^{n - r}a_{j i_{1}}  \ldots a_{i_{k} j}\, u_{1}
\ldots u_{\rho }\, a_{i_{q} i_{q + 1}}  \ldots  a_{i_{s - 1}
i_{s}}
   a_{i_{s} i_{s + 1} } \ldots a_{i_{q - 1} i_{q}}  v_{1}
\ldots v_{p}.
\]
\noindent Here by $u_{\tau}   $ and $v_{t} $ we denote  products
of coefficients whose indices form some disjoint cycles
 $(\tau = \overline {1,\rho},\;t = \overline
{1,p},\;\rho + p = r - 2)$ or there are no such products. Now for
$d$ there are the following three monomials   of ${\rm{rdet}}
_{{j}}\, {\rm {\bf A}}_{{j\,.}} ( \rm {\bf a}_{{i\,.}} )$:
\[\begin{array}{c}
     \mathord{\buildrel{\lower3pt\hbox{$\scriptscriptstyle\smile$}}\over
{d}} _{1} = ( - 1)^{n - r + 1}a_{j{\kern 1pt} i_{s + 1}}  \ldots
a_{i_{q - 1} i_{q}} \, a_{i_{q} i_{q + 1}}  \ldots a_{i_{s - 1}
i_{s}} a_{i_{s} i_{1}} \ldots a_{i_{k} j} \, u_{1} \ldots u_{\rho}
v_{1} \ldots v_{p},\\
     \mathord{\buildrel{\lower3pt\hbox{$\scriptscriptstyle\smile$}}\over {d}}
_{2} = ( - 1)^{n - r + 1}a_{j{\kern 1pt} i_{s - 1}}  \ldots
a_{i_{q + 1} i_{q}}\,  a_{i_{q} i_{q - 1}}  \ldots a_{i_{s + 1}
i_{s}} a_{i_{s} i_{1}} \ldots a_{i_{k} j}\, u_{1} \ldots u_{\rho}
v_{1} \ldots v_{p},\\
 \mathord{\buildrel{\lower3pt\hbox{$\scriptscriptstyle\smile$}}\over {d}}
_{3} = ( - 1)^{n - r}a_{j{\kern 1pt} i_{1}}  \ldots a_{i_{k} j}\,
u_{1} \ldots u_{\rho} \, a_{i_{q} i_{q - 1}}  \ldots a_{i_{s + 1}
i_{s}} a_{i_{s} i_{s - 1}}  \ldots \cdot a_{i_{q + 1} i_{q}} v_{1}
\ldots v_{p}.
\end{array}
\]

 Assume that $a_{i_{s} i_{s + 1}}   \ldots \ a_{i_{q - 1}
i_{q}}= \varphi,$ $a_{i_{q} i_{q + 1}}  \ldots  a_{i_{s - 1}
i_{s}}= \phi$, $a_{j\, i_{1}} \ldots a_{i_{k} j} =x$,\,\,$
a_{i_{q} i_{q + 1}}  \ldots  a_{i_{s - 1} i_{s}} a_{i_{s} i_{s +
1} } \ldots  a_{i_{q - 1} i_{q}}  = y$,   $ a_{i_{s} i_{s + 1}}
\ldots a_{i_{q - 1} i_{q}} a_{i_{q} i_{q + 1}} \ldots a_{i_{s - 1}
i_{s}} =y_{1}$.
 Then we
obtain $y = \phi  \varphi $,\,  $y_{1} = \varphi \phi $,
$\overline {y} = a_{i_{q} i_{q - 1}}  \ldots a_{i_{s + 1} i_{s}}
a_{i_{s} i_{s - 1}}  \ldots a_{i_{q + 1} i_{q}}  $,\,\,and
$\overline {y_{1}} = a_{i_{s} i_{s - 1}} \ldots a_{i_{q + 1}
i_{q}} a_{i_{q} i_{q - 1}} \ldots a_{i_{s + 1} i_{s}}$. Accounting
for  $a_{j\, i_{1} } = a_{i_{s} i_{1}} $, $a_{ji_{s - 1}} =
a_{i_{s} i_{s - 1}}$, $a_{ji_{s + 1}} = a_{i_{s} i_{s + 1}}  $, we
have
\[
\begin{array}{c}
\mathord{\buildrel{\lower3pt\hbox{$\scriptscriptstyle\smile$}}\over
{d}} +
\mathord{\buildrel{\lower3pt\hbox{$\scriptscriptstyle\smile$}}\over
{d}} _{1} +
\mathord{\buildrel{\lower3pt\hbox{$\scriptscriptstyle\smile$}}\over
{d}} _{2} +
\mathord{\buildrel{\lower3pt\hbox{$\scriptscriptstyle\smile$}}\over
{d}} _{3} =\\
 =( - 1)^{n - r}(xu_{1}  \ldots  u_{\rho} y
- y_{1} xu_{1}  \ldots  u_{\rho}  - \overline {y_{1}}\, x
  u_{1}  \ldots  u_{\rho}+
  xu_{1}  \ldots  u_{\rho}  \overline {y} )\times\\ \times v_{1}
\ldots  v_{p}= ( - 1)^{n - r}( xu_{1}  \ldots  u_{\rho}  t(y) -
t(y_{1} ) x
 u_{1} \ldots
u_{\rho})v_{1}  \ldots  v_{p}  = \\
 =  ( - 1)^{n - r}(t(\phi \cdot \varphi ) - t(\varphi \cdot
\phi ))xu_{1}  \ldots u_{\rho} v_{1}  \ldots v_{p}.
\end{array}
\]
Since by the rearrangement property of the trace $t(\phi \cdot
\varphi ) = t(\varphi \cdot \phi )$, then we obtain
$\mathord{\buildrel{\lower3pt\hbox{$\scriptscriptstyle\smile$}}\over
{d}} +
\mathord{\buildrel{\lower3pt\hbox{$\scriptscriptstyle\smile$}}\over
{d}} _{1} +
\mathord{\buildrel{\lower3pt\hbox{$\scriptscriptstyle\smile$}}\over
{d}} _{2} +
\mathord{\buildrel{\lower3pt\hbox{$\scriptscriptstyle\smile$}}\over
{d}} _{3} = 0$.

(iii) If the indices  $i_{s} $ and $j$ are placed  in the same
cycle, then we have the following monomials: $d_{1} ,\;\tilde
{d}_{1}
,\;\mathord{\buildrel{\lower3pt\hbox{$\scriptscriptstyle\frown$}}\over
{d}} _{1}$ or
$\mathord{\buildrel{\lower3pt\hbox{$\scriptscriptstyle\smile$}}\over
{d}} _{1} $. As shown above,   for each of them  there are another
one or three  monomials of ${\rm{rdet}} _{{j}}\, {\rm {\bf
A}}_{{j\,.}} ( {\rm {\bf a}}_{{i{\kern 1pt}.}} )$ such that the
sums of these two  or  four corresponding monomials are equal to
zero.

We have considered all possible kinds of disposition of an element
of the $i_{s} $th row as a factor of some monomial $d$ of
${\rm{rdet}} _{{j}}\, {\rm {\bf A}}_{{j\,.}} ( {\rm {\bf
a}}_{{i{\kern 1pt}.}})$. In each case there exist  one or three
corresponding monomials  for $d$ such that the sum of  two or four
monomials is equal to zero respectively. Hence, ${\rm{rdet}}
_{{j}}\, {\rm {\bf A}}_{{j\,.}} ( {\rm {\bf a}}_{{i{\kern 1pt}.}}
) = 0.$ $\blacksquare$
 \newtheorem{ns}{~~~~~\bf Corollary}[section]
\begin{ns}\label{kyrc7}
If a Hermitian matrix ${\rm {\bf A}}\in {\rm M}\left( {n,\bf H}
\right)$ consists two same rows (columns), then $\det{\rm {\bf
A}}=0$.
\end{ns}
{\textit{Proof.}} Suppose the $i$th row of ${\rm {\bf A}}$
coincides with the $j$th row, i.e. $a_{ik} = a_{jk} $\, $\forall k
\in I_{n} $, ${\left\{ {i,j} \right\}} \in I_{n} $, $i \ne j$.
Then $\overline {a_{ik}}  = \overline {a_{jk}}  $, $(\forall k \in
I_{n})$. Since the matrix ${\rm {\bf A}}$ is Hermitian, we get
 $\forall k \in
I_{n}$ $a_{ki} = a_{kj}$, where ${\left\{ {i,j} \right\}} \in
I_{n} $, $i \ne j$. This means that if a Hermitian matrix has two
same rows, then it has two same corresponding columns as well. The
matrix ${\rm {\bf A}}$ may be represented as ${\rm {\bf A}}_{j.}
\left( {{\rm {\bf a}}_{i.} } \right)$, where the matrix ${\rm {\bf
A}}_{j.} \left( {{\rm {\bf a}}_{i.} } \right)$ is obtained from
${\rm {\bf A}}$ by replacing the $j$th row with the $i$th row.
Then by Theorem \ref{kyrc6}, we have
\[
\det {\rm {\bf A}} = {\rm{rdet}} _{i} {\rm {\bf A}} = {\rm{rdet}}
_{i} {\rm {\bf A}}_{j.} \left( {{\rm {\bf a}}_{i.}}  \right) =
0.\blacksquare
\]
\begin{theorem} \label{kyrc8}If the matrix ${\rm {\bf A}}_{.{\kern 1pt} i}
\left( {{\rm {\bf a}}_{.j}}  \right)$ is obtained from a Hermitian
matrix ${\rm {\bf A}}\in {\rm M}\left( {n,\bf H} \right)$ by
replacing of its $i$th column with the $j$th column, then
${\rm{cdet}} _{{i}}\, {\rm {\bf A}}_{{.\,i}} ( {\rm {\bf
a}}_{{.j}}) = 0$, $(\forall i,j=\overline{1,n},\,i\neq j).$
\end{theorem}
{\textit{Proof.}} The proof of this theorem is analogous to that
of Theorem \ref{kyrc6}.

 From Theorems \ref{kyrc6}, \ref{kyrc8} and basic
properties of  row and column determinants for arbitrary matrices
we have the following theorems.
\begin{theorem}\label{kyrc22} If the matrix ${\rm {\bf
A}}_{i{\kern 1pt}.} \left( {b \cdot {\rm {\bf a}}_{j{\kern 1pt}.}}
\right)$ is obtained from a Hermitian matrix ${\rm {\bf A}}\in
M\left( {n,\bf H} \right)$   by replacing of its $i$th row with
the $j$th row multiplied by $b \in \bf H$ on the left, then
${\rm{rdet}}_{i}\, {\rm {\bf A}}_{i{\kern 1pt}.} \left( {b \cdot
{\rm {\bf a}}_{j{\kern 1pt}.}} \right) = 0$, $(\forall
i,j=\overline{1,n},\,i\neq j).$
\end{theorem}
\begin{theorem} If the matrix ${\rm {\bf
A}}_{.{\kern 1pt}j} \left( {{\rm {\bf a}}_{.{\kern 1pt} i} \cdot
b} \right)$ is obtained from a Hermitian matrix ${\rm {\bf A}}\in
{\rm M}\left( {n,\bf H} \right)$   by replacing of its $j$th
column with the $i$th column right-multiplied by $b \in\bf H$,
then ${\rm{cdet}} _{j}\, {\rm {\bf A}}_{.{\kern 1pt}j} \left(
{{\rm {\bf a}}_{.{\kern 1pt} i} \cdot b} \right) = 0$, $(\forall
i,j=\overline{1,n},\,i\neq j).$
\end{theorem}
\begin{theorem} \label{kyrch23}
If the matrix ${\rm {\bf A}}_{.j} \left( {{\rm {\bf a}}_{.{\kern
1pt} {\kern 1pt} i} \cdot b} \right)$ is obtained from a Hermitian
matrix ${\rm {\bf A}}\in {\rm M}\left( {n,\bf H} \right)$ by
replacing of its $j$th column with the $i$th column multiplied by
$b \in \bf H$ on the right, then ${\rm{rdet}}_{j} \,{\rm {\bf
A}}_{.{\kern 1pt}j} \left( {{\rm {\bf a}}_{.{\kern 1pt} i} \cdot
b} \right) = 0$, $(\forall i,j=\overline{1,n},\,i\neq j).$
\end{theorem}
{\textit{Proof.}} We assume $n > 3$ for  ${\rm {\bf A}}\in {\rm
M}\left( {n,\bf H}\right)$. The case $n \le 3$ can be  easily
proved by direct calculation. Consider some monomial $d$ of
${\rm{rdet}} _{{j}}\, {\rm {\bf A}}_{.j} \left( {{\rm {\bf
a}}_{.\,i} \cdot b} \right)$. Suppose indices of its coefficients
form a permutation as a product of $r$ disjoint cycles, and denote
$i = i_{s} $. Consider all
 possibilities of disposition of an entry of the $i_{s} $th
row  in the monomial $d$.

\noindent (i) Suppose an entry of the $i_{s} $th row is placed in
$d$ such that the index $i_{s} $ opens some disjoint cycle, i.e.:
\begin{equation}
\label{eq1} d = ( - 1)^{n - r}a_{j{\kern 1pt} i_{1}}  \ldots
a_{i_{k} j} \,b{\kern 1pt} {\kern 1pt} u_{1} \ldots u_{\rho}
\,a_{i_{s} i_{s + 1}}  \ldots a_{i_{s + m} i_{s}}  v_{1} \ldots
v_{p} {\rm ,}
\end{equation}

\noindent Here we denote by $u_{\tau}  $ and $v_{t} $  products of
coefficients whose indices  form  disjoint cycles $(\forall \tau =
\overline {1,\rho} ,\;\forall t = \overline {1,p},\, \rho + p = r
- 2)$ or there are no such products. For $d$ there are the
following three monomials  of ${\rm{rdet}} _{{j}}\, {\rm {\bf
A}}_{.j} \left( {{\rm {\bf a}}_{.\,i} \cdot b} \right)$.
\[\begin{array}{c}
d_{1} = ( - 1)^{n - r}a_{j{\kern 1pt} i_{1}}  \cdot \ldots \cdot
a_{i_{k} j} \cdot b \cdot u_{1} \cdot \ldots \cdot u_{\rho}  \cdot
a_{i_{s} i_{s + m}} \cdot \ldots \cdot a_{i_{s + 1} i_{s}}  \cdot
v_{1} \cdot \ldots \cdot v_{p},
\\
d_{2} = ( - 1)^{n - r + 1}a_{j{\kern 1pt} i_{1}}  \cdot \ldots
\cdot a_{i_{k} i_{s}}  \cdot a_{i_{s} i_{s + 1}}  \ldots a_{i_{s +
m} j} \cdot b \cdot u_{1} \cdot \ldots \cdot u_{\rho}  \cdot v_{1}
\cdot \ldots \cdot v_{p},
\\
d_{3} = ( - 1)^{n - r + 1}a_{ji_{1}}  \cdot \ldots \cdot a_{i_{k}
i_{s}} \cdot a_{i_{s} {\kern 1pt} i_{s + m}}  \cdot \ldots \cdot
a_{i_{s + 1} j} \cdot b \cdot u_{1} \cdot \ldots \cdot u_{\rho}
\cdot v_{1} \cdot \ldots \cdot v_{p}.
\end{array}
\]
If $a_{j{\kern 1pt} i_{1}}  \cdot \ldots \cdot a_{i_{k} j} = x$,
$a_{i_{s} i_{s + 1}}  \ldots a_{i_{s + m} i_{s}}  = y$, then
$\overline {y} = a_{i_{s} i_{s + m}} \ldots a_{i_{s + 1} i_{s} }
$. Taking into account  $a_{i_{k} j{\kern 1pt}}  = a_{i_{k} i_{s}}
$, $a_{i_{s + m} j} = a_{i_{s + m} i_{s}}$, $a_{i_{s + 1} j} =
a_{i_{s + 1} i_{s}}  $, we consider the sum of these monomials.
\[\begin{array}{c}
d + d_{1} + d_{2} + d_{3} =
\\
 = {\rm (} - 1{\rm )}^{n - r}{\rm (}x \cdot b \cdot u_{1} \cdot \ldots \cdot
u_{\rho}  \cdot y + x \cdot b \cdot u_{1} \cdot \ldots \cdot
u_{\rho}  \cdot \overline {y} - x \cdot y \cdot b \cdot u_{1}
\cdot \ldots \cdot u_{\rho}  -
\\
 - x \cdot \overline {y} \cdot b \cdot u_{1} \cdot \ldots \cdot u_{\rho}
{\rm )} \cdot v_{1} \cdot \ldots \cdot v_{p} = {\rm (} - 1{\rm
)}^{n - r}{\rm (}x \cdot b \cdot u_{1} \cdot \ldots \cdot u_{\rho}
\cdot {\rm (}y + \overline {y} {\rm )} -
\\
 - x \cdot {\rm (}y + \overline {y} {\rm )} \cdot b \cdot u_{1} \cdot \ldots
\cdot u_{\rho}  {\rm )} \cdot v_{1} \cdot \ldots \cdot v_{p} =
{\rm (} - 1{\rm )}^{n - r}{\rm (}x \cdot b \cdot u_{1} \cdot
\ldots \cdot u_{\rho} \cdot t{\rm (}y{\rm )} -
\\
 - x \cdot t{\rm (}y{\rm )} \cdot b \cdot u_{1} \cdot \ldots \cdot u_{\rho}
{\rm )} \cdot v_{1} \cdot \ldots \cdot v_{p} = 0{\rm .}
\end{array}
\]
Thus among the monomials of  ${\rm {rdet}} _{j} {\rm {\bf A}}_{.j}
\left( {{\rm {\bf a}}_{.\,i} \cdot b} \right)$ we find
 three monomials  for $d$ such  that the sum of these  monomials and
 $d$ is equal to zero.
If in (\ref{eq1}) $m = 0$ or $m = 1$, we  obtain such monomials
accordingly:
\[\begin{array}{c}
\tilde {d} = ( - 1)^{n - r}a_{j{\kern 1pt} i_{1}}  \cdot \ldots
\cdot a_{i_{k} j} \cdot b \cdot u_{1} \cdot \ldots \cdot u_{\rho}
\cdot a_{i_{s} i_{s}}  \cdot v_{1} \cdot \ldots \cdot v_{p},
\\
\mathord{\buildrel{\lower3pt\hbox{$\scriptscriptstyle\smile$}}\over
{d}} = ( - 1)^{n - r}a_{j{\kern 1pt} i_{1}}  \cdot \ldots \cdot
a_{i_{k} j} \cdot b \cdot u_{1} \cdot \ldots \cdot u_{\rho}  \cdot
a_{i_{s} i_{s + 1}}  \cdot a_{i_{s + 1} i_{s}}  \cdot v_{1} \cdot
\ldots \cdot v_{p}.
\end{array}
\]
There are the following  monomials for them.
\[\begin{array}{c}
\tilde {d}_{1} = ( - 1)^{n - r + 1}a_{j{\kern 1pt} i_{1}}  \cdot
\ldots \cdot a_{i_{k} i_{s}}  a_{i_{s} j} \cdot b \cdot u_{1}
\cdot \ldots \cdot u_{\rho}  \cdot v_{1} \cdot \ldots \cdot v_{p},
\\
\mathord{\buildrel{\lower3pt\hbox{$\scriptscriptstyle\smile$}}\over
{d}} _{1} = ( - 1)^{n - r + 1}a_{j{\kern 1pt} i_{1}}  \cdot \ldots
\cdot a_{i_{k} i_{s}}  \cdot a_{i_{s} i_{s + 1}}  \cdot a_{i_{s +
1} j} \cdot b \cdot u_{1} \cdot \ldots \cdot u_{\rho}  \cdot v_{1}
\cdot \ldots \cdot v_{p},
\end{array}
\]
Taking into account $a_{i_{k} j{\kern 1pt}}  = a_{i_{k} i_{s}}  $,
$a_{i_{s} j} = a_{i_{s} i_{s}}  {\rm ,}_{{\rm} }a_{i_{s + 1} j} =
a_{i_{s + 1} i_{s}} $, and $a_{i_{s} i_{s}}  \in {\bf R}$,
$a_{i_{s} i_{s + 1}}  a_{i_{s + 1} i_{s}}  = n\left( {a_{i_{s}
i_{s + 1}} } \right) \in {\bf R}$, we get
\[\begin{array}{c}
\tilde {d} + \tilde {d}_{1} = ( - 1)^{n - r}{\rm (}a_{j{\kern 1pt}
i_{1}} \cdot \ldots \cdot a_{i_{k} j} \cdot b \cdot u_{1} \cdot
\ldots \cdot u_{\rho}  \cdot a_{i_{s} i_{s}}  -
\\
 - a_{j{\kern 1pt} i_{1}}  \cdot \ldots \cdot a_{i_{k} i_{s}}  \cdot
a_{i_{s} j} \cdot b \cdot u_{1} \cdot \ldots \cdot u_{\rho}  {\rm
)} \cdot v_{1} \cdot \ldots \cdot v_{p} = 0{\rm ,}
\end{array}
\]
\[\begin{array}{c}
\mathord{\buildrel{\lower3pt\hbox{$\scriptscriptstyle\smile$}}\over
{d}} +
\mathord{\buildrel{\lower3pt\hbox{$\scriptscriptstyle\smile$}}\over
{d}} _{1} = ( - 1)^{n - r}{\rm (}a_{j{\kern 1pt} i_{1}}  \cdot
\ldots \cdot a_{i_{k} j} \cdot b \cdot u_{1} \cdot \ldots \cdot
u_{\rho}  \cdot n\left( {a_{i_{s} i_{s + 1}} }  \right) -
\\
 - a_{j{\kern 1pt} i_{1}}  \cdot \ldots \cdot a_{i_{k} i_{s}}  \cdot n\left(
{a_{i_{s} i_{s + 1}} }  \right) \cdot b \cdot u_{1} \cdot \ldots
\cdot u_{\rho}  {\rm )} \cdot v_{1} \cdot \ldots \cdot v_{p} =
{\rm 0}{\rm .}
\end{array}
\]
Hence,  the sums of corresponding two monomials of ${\rm{rdet}}
_{j} {\rm {\bf A}}_{.j} \left( {{\rm {\bf a}}_{.\,i} \cdot b}
\right)$ are equal to zero in this case.

\noindent
 (ii) Now suppose that the index $i_{s} $  is placed in  another disjoint
  cycle  than the index $j$ and does not open this cycle,
\[
  \hat {d}= ( - 1)^{n - r}a_{j{\kern 1pt} i_{1}}  \ldots a_{i_{k} j}
b\,u_{1} \ldots u_{\rho}  a_{i_{q} i_{q + 1}}  \ldots a_{i_{s - 1}
i_{s}}  a_{i_{s} i_{s + 1}}  \ldots a_{i_{q - 1} i_{q}}  v_{1}
\ldots v_{p}.
\]
\noindent Here we denote by $u_{\tau}   $ and $v_{t} $  products
of coefficients whose indices   form  disjoint cycles
 $(\tau = \overline {1,\rho},\;t = \overline
{1,p},\;\rho + p = r - 2)$ or there are no such products. Now for
$d$ we have the following three monomials   of ${\rm{rdet}} {\rm
{\bf A}}_{.j} \left( {{\rm {\bf a}}_{.\,i} \cdot b} \right)$.
\[\begin{array}{c}
\hat {d}_{1} = ( - 1)^{n - r}a_{j{\kern 1pt} i_{1}}  \ldots
a_{i_{k} j} bu_{1} \ldots u_{\rho}  a_{i_{q} i_{q - 1}}  \ldots
a_{i_{s + 1} i_{s}} a_{i_{s} i_{s - 1}}  \ldots a_{i_{q + 1}
i_{q}}  v_{1} \ldots v_{p},
\\
\hat {d}_{2} = ( - 1)^{n - r}a_{j{\kern 1pt} i_{1}}  \ldots
a_{i_{k} i_{s}} a_{i_{s} i_{s - 1}}  \ldots \cdot a_{i_{q + 1}
i_{q}}  a_{i_{q} i_{q - 1}} \ldots a_{i_{s + 1} j} \,bu_{1} \ldots
u_{\rho}  v_{1} \ldots v_{p},
\\
\hat {d}_{3} = ( - 1)^{n - r}a_{j{\kern 1pt} i_{1}}  \ldots
a_{i_{k} i_{s}} a_{i_{s} i_{s - 1}}  \ldots a_{i_{q + 1} i_{q}}
a_{i_{q} i_{q - 1}}  \ldots a_{i_{s + 1} j} \,b\,u_{1} \ldots
u_{\rho}  v_{1} \ldots v_{p} \end{array}
\]
 Assume that $a_{j{\kern 1pt} i_{1}}  \cdot \ldots \cdot a_{i_{k} j} = x$,
 $a_{i_{q} i_{q + 1}}  \cdot \ldots \cdot a_{i_{s - 1} i_{s}}
= \phi $, $a_{i_{s} i_{s + 1}}  \cdot \ldots \cdot a_{i_{q - 1}
i_{q}}  = \varphi $, then $a_{i_{s} i_{s - 1}}  \cdot \ldots \cdot
a_{i_{q + 1} i_{q}}  = \overline {\phi}  \,$ è $a_{i_{q} i_{q -
1}}  \cdot \ldots \cdot a_{i_{s + 1} i_{s}}  = \overline
{\varphi}$. Taking into account  $a_{i_{k} j{\kern 1pt}} =
a_{i_{k} i_{s}}  $, $a_{i_{s - 1} j} = a_{i_{s - 1} i_{s}}
,_{{\rm} }a_{i_{s + 1} j} = a_{i_{s + 1} i_{s} } $, we have
\[\begin{array}{c}
\hat {d}+ \hat {d}_{1} + \hat {d}_{2} + \hat {d}_{3} =
\\
 = {\rm (} - 1{\rm )}^{n - r}{\rm (}xbu_{1} \ldots u_{\rho}  \phi \,\varphi
+ xbu_{1} \ldots u_{\rho}  \overline {\varphi}  \,\overline {\phi}
-
\\
 - x\varphi \,\phi bu_{1} \ldots u_{\rho}  - x\overline {\phi \,} \overline
{\varphi}  \,bu_{1} \ldots u_{\rho}  {\rm )}v_{1} \ldots v_{p} =
\\
 = {\rm (} - 1{\rm )}^{n - r}{\rm (}xbu_{1} \ldots u_{\rho}  {\rm (}\phi
\,\varphi + \overline {\phi \varphi \,} {\rm )} - x{\rm (}\varphi
\,\phi + \overline {\varphi \phi}  \,{\rm )}bu_{1} \ldots u_{\rho}
{\rm )}v_{1} \ldots v_{p} =
\\
 = {\rm (} - 1{\rm )}^{n - r}{\rm (}xbu_{1} \ldots u_{\rho}  t{\rm (}\phi
\,\varphi {\rm )} - x\,t{\rm (}\varphi \,\phi {\rm )}bu_{1} \ldots
u_{\rho} {\rm )}v_{1} \ldots v_{p}.
\end{array}\]
Since  $t(\phi \cdot \varphi ) = t(\varphi \cdot \phi )$ by the
rearrangement property of the trace, then we obtain $\hat {d} +
\hat {d}_{1} + \hat {d}_{2} + \hat {d}_{3} = 0$.

\noindent (iii) If the indices  $i_{s} $ and $j$ are placed  in
the same cycle, then we have the following monomials: $d_{2} $ or
$d_{3} $, $\hat {d}_{2} $ or $\hat {d}_{3} $, and either $\tilde
{d}_{1} $ or
$\mathord{\buildrel{\lower3pt\hbox{$\scriptscriptstyle\smile$}}\over
{d}} _{1} $. As shown above,   for each of them  there are another
one or three  monomials of ${\rm{rdet}} _{j} {\rm {\bf A}}_{.j}
\left( {{\rm {\bf a}}_{.\,i} \cdot b} \right)$ such that the sums
of these two  or  four corresponding monomials are equal to zero.

We have considered all possible kinds of disposition of an entry
of the $i_{s} $th row as a factor of some monomial $d$ of
${\rm{rdet}} _{j} {\rm {\bf A}}_{.j} \left( {{\rm {\bf a}}_{.\,i}
\cdot b} \right)$. In each case there exist  one or three
corresponding monomials  for $d$ such that the sums of these two
or four monomials are equal to zero respectively. Hence,
${\rm{rdet}}_{j} {\rm {\bf A}}_{.j} \left( {{\rm {\bf a}}_{.\,i}
\cdot b} \right) = 0$, $\left( {\forall i{\rm ,}j = \overline
{1{\rm ,}n} {\rm ,}\,\,i \ne j} \right)$. $\blacksquare$

\begin{ns} If the matrix  ${\rm {\bf A}}_{.j} \left( {{\rm {\bf a}}_{.{\kern
1pt} {\kern 1pt} i}}  \right)$  is obtained from a Hermitian
matrix ${\rm {\bf A}}\in {\rm M}\left( {n,\bf H} \right)$ by
replacing of its $j$th column with the $i$th column, then
${\rm{rdet}} _{{j}}\, {\rm {\bf A}}_{.j} \left( {{\rm {\bf
a}}_{.\,i}} \right)= 0$, $\left( {\forall i{\rm ,}j = \overline
{1{\rm ,}n} {\rm ,}\,\,i \ne j} \right)$.
\end{ns}
{\textit{Proof}}. The proof of this lemma follows from Theorem
\ref{kyrch23} by put $b = 1$.
\begin{theorem} \label{kyrch24}
If the matrix ${\rm {\bf A}}_{i\,.} \left( {b \cdot {\rm {\bf
a}}_{j.}}  \right)$ is obtained from a Hermitian matrix ${\rm {\bf
A}}\in {\rm M}\left( {n,\bf H} \right)$ by replacing of its $i$th
row with the $j$th row multiplied by $b \in \bf H$ on the left,
then ${\rm{cdet}} _{i} {\rm {\bf A}}_{i\,.} \left( {b \cdot {\rm
{\bf a}}_{j.}}  \right) = 0$, $\left( {\forall i,j = \overline
{1,n},\,\,i \ne j} \right).$
\end{theorem}
{\textit{Proof}}. The proof of this theorem is analogous to that
of Theorem \ref{kyrch23}.
\begin{ns} If the matrix  ${\rm {\bf A}}_{i\,.} \left( {{\rm {\bf a}}_{j.}
} \right)$ is obtained from a Hermitian matrix ${\rm {\bf A}}\in
{\rm M}\left( {n,\bf H} \right)$ by replacing of its $i$th row
with the $j$th, then  ${\rm{cdet}} _{i} {\rm {\bf A}}_{i\,.}
\left( {{\rm {\bf a}}_{j.}}  \right) = 0$, $\left( {\forall i,j =
\overline {1,n},\,\,i \ne j} \right).$
\end{ns}
{\textit{Proof}}. The proof of this corollary follows from Theorem
\ref{kyrch24} by put $b = 1$.

 The following theorems immediately follows from the previous theorems
 and basic
properties of the row and column determinants for arbitrary
matrices.
\begin{theorem}\label{kyrch25} If the $i$th row of
a Hermitian matrix ${\rm {\bf A}}\in {\rm M}\left( {n,\bf H}
\right)$ is replaced with a left linear combination of its other
rows, i.e. ${\rm {\bf a}}_{i.} = c_{1} {\rm {\bf a}}_{i_{1} .} +
\ldots + c_{k}  {\rm {\bf a}}_{i_{k} .}$, where $ c_{l} \in {\bf
H}$ for $\forall l = \overline {1,k}$ and $\{i,i_{l}\}\subset
I_{n} $, then
\[
 {\rm{rdet}}_{i}\, {\rm {\bf A}}_{i \, .} \left(
{c_{1} {\rm {\bf a}}_{i_{1} .} + \ldots + c_{k} {\rm {\bf
a}}_{i_{k} .}}  \right) = {\rm{cdet}} _{i}\, {\rm {\bf A}}_{i\, .}
\left( {c_{1}
 {\rm {\bf a}}_{i_{1} .} + \ldots + c_{k} {\rm {\bf
a}}_{i_{k} .}}  \right) = 0.
\]
\end{theorem}
\begin{ns} \label{kyrch26}
If some  row of a Hermitian matrix ${\rm {\bf A}}\in {\rm M}\left(
{n,\bf H} \right)$ is a left linear combination of its other rows,
then $ det{\rm {\bf A}}=0$.
\end{ns}
{\textit{Proof}} Let ${\rm {\bf a}}_{i.} = c_{1} {\rm {\bf
a}}_{i_{1} .} + \ldots + c_{k}  {\rm {\bf a}}_{i_{k} .}$, where $
c_{l} \in {\bf H}$ for $\forall l = \overline {1,k}$ and
$\{i,i_{l}\}\subset I_{n} $. The matrix ${\rm {\bf A}}$ may be
represented as ${\rm {\bf A}}_{i.} \left( c_{1} {\rm {\bf
a}}_{i_{1} .} + \ldots + c_{k}  {\rm {\bf a}}_{i_{k} .} \right)$.
Then by Theorem \ref{kyrch25}, we have
\[
\det {\rm {\bf A}} = {\rm{rdet}} _{i} {\rm {\bf A}} = {\rm{rdet}}
_{i} {\rm {\bf A}}_{i.} \left( c_{1} {\rm {\bf a}}_{i_{1} .} +
\ldots + c_{k}  {\rm {\bf a}}_{i_{k} .} \right) = 0. \blacksquare
\]
\begin{theorem}\label{kyrch27} If the $j$th column of
 a Hermitian matrix ${\rm {\bf A}}\in
{\rm M}\left( {n,\bf H} \right)$   is replaced with a right linear
combination of its other columns, i.e. ${\rm {\bf a}}_{.j} = {\rm
{\bf a}}_{.j_{1}}   c_{1} + \ldots + {\rm {\bf a}}_{.j_{k}} c_{k}
$, where $c_{l} \in{\bf H}$ for $\forall l = \overline {1,k}$ and
$\{j,j_{l}\}\subset J_{n}$,  then
 \[{\rm{cdet}} _{j}\, {\rm {\bf A}}_{.j}
\left( {{\rm {\bf a}}_{.j_{1}} c_{1} + \ldots + {\rm {\bf
a}}_{.j_{k}}c_{k}} \right) ={\rm{rdet}} _{j} \,{\rm {\bf A}}_{.j}
\left( {{\rm {\bf a}}_{.j_{1}}  c_{1} + \ldots + {\rm {\bf
a}}_{.j_{k}}  c_{k}} \right) = 0.
\]
\end{theorem}
\begin{ns}\label{kyrch28}
If some  column of a Hermitian matrix ${\rm {\bf A}}\in {\rm
M}\left( {n,\bf H} \right)$ is a right linear combination of its
other columns, then $ det{\rm {\bf A}}=0$.
\end{ns}
{\textit{Proof}}. The proof of this corollary is analogous to that
of Corollary \ref{kyrch26} and follows from Theorem \ref{kyrch27}.
\begin{theorem} If the $i$th row of a Hermitian matrix ${\rm {\bf A}}\in
{\rm M}\left( {n,\bf H} \right)$  is added a left linear
combination of its other rows, then \[\begin{array}{c}
   {\rm{rdet}}_{i} \,{\rm {\bf A}}_{i \cdot } \left( {{\rm {\bf
a}}_{i.} + c_{1} \cdot {\rm {\bf a}}_{i_{1} .} + \ldots + c_{k}
\cdot {\rm {\bf a}}_{i_{k} .}} \right) =\\
  ={\rm{cdet}} _{i}
\,{\rm {\bf A}}_{i \cdot} \left( {{\rm {\bf a}}_{i.} + c_{1} \cdot
{\rm {\bf a}}_{i_{1} .} + \ldots + c_{k} \cdot {\rm {\bf
a}}_{i_{k} .}}  \right) = \det {\rm {\bf A}},
\end{array}\]
 where $ c_{l} \in {\bf H}$ for $\forall l = \overline {1,k}$ and
$\{i,i_{l}\}\subset I_{n}$.
\end{theorem}
\begin{theorem} If the $j$th column of a Hermitian matrix ${\rm {\bf A}}\in
{\rm M}\left( {n,\bf H} \right)$   is added a right linear
combination of its other columns, then
\[\begin{array}{c}
   {\rm{cdet}} _{j}\, {\rm {\bf A}}_{.j} \left( {{\rm {\bf a}}_{.j} +
{\rm {\bf a}}_{.j_{1}}
 c_{1} + \ldots + {\rm {\bf a}}_{.j_{k}}   c_{k}}
\right) = \\
  ={\rm{rdet}} _{j} \,{\rm {\bf A}}_{.j} \left( {{\rm {\bf
a}}_{.j} + {\rm {\bf a}}_{.j_{1}}   c_{1} + \ldots + {\rm {\bf
a}}_{.j_{k}}  c_{k}}  \right) = \det {\rm {\bf A}},
\end{array}\]
where $c_{l} \in{\bf H}$ for $\forall l = \overline {1,k}$ and
$\{j,j_{l}\}\subset J_{n}$.
\end{theorem}

\section{ The inverse of a Hermitian matrix}
\begin{definition} A Hermitian matrix ${\rm {\bf A}}\in
{\rm M}\left( {n,\bf H} \right)$ is called nonsingular if $\det
{\rm {\bf A}} \ne 0$.
\end{definition}
\begin{theorem} \label{kyrc9} There exists a unique right
inverse  matrix $(R{\rm {\bf A}})^{ - 1}$ and a unique left
inverse matrix $(L{\rm {\bf A}})^{ - 1}$ of a nonsingular
Hermitian matrix ${\rm {\bf A}}\in {\rm M}\left( {n,\bf H}
\right)$ such that $\left( {R{\rm {\bf A}}} \right)^{ - 1} =
\left( {L{\rm {\bf A}}} \right)^{ - 1} = :{\rm {\bf A}}^{ - 1}$,
where
\[\begin{array}{c}
   \left( {R{\rm {\bf A}}} \right)^{ - 1} = {\frac{{1}}{{\det
{\rm {\bf A}}}}}
\begin{pmatrix}
  R_{11} & R_{21} & \cdots & R_{n1}\\
  R_{12} & R_{22} & \cdots & R_{n2}\\
  \cdots & \cdots & \cdots& \cdots\\
  R_{1n} & R_{2n} & \cdots & R_{nn}
\end{pmatrix},\\
  \left( {L{\rm {\bf A}}} \right)^{ - 1} = {\frac{{1}}{{\det {\rm
{\bf A}}}}}
\begin{pmatrix}
  L_{11} & L_{21} & \cdots & L_{n1} \\
  L_{12} & L_{22} & \cdots & L_{n2} \\
  \cdots & \cdots & \cdots & \cdots \\
  L_{1n} & L_{2n} & \cdots & L_{nn}
\end{pmatrix},
\end{array}
\]
 $R_{ij}$,  $L_{ij}$ are right and left $ij$-th cofactor of ${\rm {\bf
 A}}$ respectively, $\left( {\forall i,j =
\overline {1,n}} \right)$.
\end{theorem}
 {\textit{Proof.}} Let ${\rm {\bf B}} = {\rm {\bf A}} \cdot \left( {R{\rm
{\bf A}}} \right)^{ - 1}$.  We obtain the entries of ${\rm {\bf
B}}$ by direct calculations.
\[\begin{array}{c}
   b_{i{\kern 1pt} i} = \left( {\det {\rm {\bf A}}} \right)^{ -
1}{\sum\limits_{j = 1}^{n} {a_{i{\kern 1pt} j} \cdot R_{i{\kern
1pt} j}} } = \left( {\det {\rm {\bf A}}} \right)^{ -
1}{\rm{rdet}}_{i}\, {\rm {\bf A}} = {\frac{{\det {\rm {\bf
A}}}}{{\det {\rm {\bf A}}}}} = 1{\rm ,} \left( {\forall i =
\overline {1,n}} \right),\\
  b_{i{\kern 1pt} j} = \left( {\det {\rm {\bf A}}} \right)^{ -
1}{\sum\limits_{s = 1}^{n} {a_{i{\kern 1pt} s} \cdot R_{j{\kern
1pt} s}} } = \left( {\det {\rm {\bf A}}} \right)^{ - 1}{\rm{rdet}}
_{j} {\rm {\bf A}}_{j{\kern 1pt}.} \left( {{\rm {\bf a}}_{i{\kern
1pt}.}} \right), \quad \left( {i \ne j} \right).
\end{array}
\]
If $i \ne j$, then by Theorem \ref{kyrc6} ${\rm{rdet}} _{j} {\rm
{\bf A}}_{j{\kern 1pt}.} \left( {{\rm {\bf a}}_{i{\kern 1pt}.}}
\right)=0$. Consequently $b_{i{\kern 1pt} j} = 0$. Thus ${\rm {\bf
B}} = {\rm {\bf I}}$ and $\left( {R{\rm {\bf A}}} \right)^{ - 1}$
is the right inverse of the Hermitian matrix ${\rm {\bf A}}$.

Suppose now that ${\rm {\bf D}} = \left( {L{\rm {\bf A}}}
\right)^{ - 1}{\rm {\bf A}}$. We again  get the entries of ${\rm
{\bf D}}$ by multiplying matrices.
\[\begin{array}{c}
 d_{i{\kern 1pt} i} = \left( {\det {\rm {\bf A}}} \right)^{ -
1}{\sum\limits_{i = 1}^{n} {L_{i{\kern 1pt} j} \cdot a_{i{\kern
1pt} j}} } = \left( {\det {\rm {\bf A}}} \right)^{ - 1}{\rm{cdet}}
_{j} {\rm {\bf A}} = {\frac{{\det {\rm {\bf A}}}}{{\det {\rm {\bf
A}}}}} = 1, \left( {\forall i = \overline {1,n}}  \right),\\
  d_{i{\kern 1pt} j} = \left( {\det {\rm {\bf A}}} \right)^{ -
1}{\sum\limits_{s = 1}^{n} {L_{s{\kern 1pt} {\kern 1pt} i} \cdot
a_{s{\kern 1pt} j}} }  = \left( {\det {\rm {\bf A}}} \right)^{ -
1}{\rm{cdet}} _{i} {\rm {\bf A}}_{i.} \left( {{\rm {\bf a}}_{j.}}
\right), \left( {i \ne j} \right).
\end{array}
\]
If $i \ne j$, then by Theorem \ref{kyrc8} ${\rm{cdet}} _{i} {\rm
{\bf A}}_{i.} \left( {{\rm {\bf a}}_{j.}} \right)=0$. Therefore
$d_{i{\kern 1pt} j} = 0$. Thus ${\rm {\bf D}} = {\rm {\bf I}}$ and
$\left( {L{\rm {\bf A}}} \right)^{ - 1}$ is the left inverse
 of the Hermitian matrix ${\rm {\bf A}}$.

The equality $\left( {R{\rm {\bf A}}} \right)^{ - 1} = \left(
{L{\rm {\bf A}}} \right)^{ - 1}$ is immediate from the well-known
fact that if there exists an inverse matrix over an arbitrary skew
field, then it is unique.$\blacksquare$

\section{Properties of the corresponding\\ Hermitian matrices}
 Denote by ${\bf H}^{m\times n}$  the set of $m\times n$ matrices
  with entries in ${\bf H}$.
  \begin{definition}
If ${\rm {\bf A}}\in {\bf H}^{m\times n} $, then  the matrix ${\rm
{\bf A}}^{ * }{\rm {\bf A}} \in {\rm M}\left( {n,\bf H} \right) $
is called its  left corresponding Hermitian and ${\rm {\bf A}}{\rm
{\bf A}}^{ *} \in {\rm M}\left( {m,\bf H} \right)$ is called its
right corresponding Hermitian matrix.
\end{definition}
\begin{theorem}  \label{kyrc10} If an arbitrary
column of ${\rm {\bf A}}\in {\bf H}^{m\times n} $  is a right
linear combination of its other columns, then $\det {\rm {\bf
A}}^{ * }{\rm {\bf A}} = 0.$
\end{theorem}
 {\textit{Proof.}} Let the $j$th column of ${\rm {\bf A}}\in {\bf H}^{m\times n} $  be a right
linear combination of its other columns. That is ${{\rm {\bf
a}}_{.j} = {\rm {\bf a}}_{.j_{1}} c_{1} + \ldots + {\rm {\bf
a}}_{.j_{k}} c_{k}}$, where $c_{l} \in{\bf H}$ for $\forall l =
\overline {1,k}$ and $\{j,j_{l}\}\subset J_{n}$. Then the $j$th
row of ${\rm {\bf A}}^{ *}$ is the left linear combination of its
rows, ${{\rm {\bf a}}_{.j}^{ *} = \overline{c_{1}}\,{\rm {\bf
a}}_{.j_{1}}^{ *}  + \ldots + \overline{c_{k}}\,{\rm {\bf
a}}_{.j_{k}}^{ *} }$. Consider the Hermitian matrix $ {\rm {\bf
A}}^{ * }{\rm {\bf A}}$. It is easy to verify that the $j$th
column of $ {\rm {\bf A}}^{ * }{\rm {\bf A}}$ is a right linear
combination of its other columns. Therefore by Corollary
\ref{kyrch28} we have $\det {\rm {\bf A}}^{ * }{\rm {\bf
A}}=0.\blacksquare$
\begin{theorem}\label{kyrc11} If some row of ${\rm {\bf
A}}\in {\bf H}^{m\times n} $  is a left linear combination of its
other rows, then $\det {\rm {\bf A}}{\rm {\bf A}}^{ *} = 0.$
\end{theorem}
{\textit{Proof}}. The proof of this theorem is analogous to that
of Theorem \ref{kyrc10} and follows from Corollary \ref{kyrch26}.
\begin{remark} Since the principal submatrices of a Hermitian matrix over ${\bf H}$
 are
Hermitian, then the basis principal minor may be defined  in this
noncommutative case as well.
\end{remark}
\begin{definition}
The basis principal minor of a Hermitian matrix over ${\bf H}$ is
defined as the nonzero determinant of the largest principal
submatrix. Then  rows and  columns included in the basis principal
minor are called the basis ones as well.
\end{definition}
\begin{definition} If  rows and  columns with indices $i_{1} ,\ldots
,i_{r} $ of  ${\rm {\bf A}}^{ *} {\rm {\bf A}}$ are basis, then
 rows with indices $i_{1} ,\ldots ,i_{r} $ of  ${\rm {\bf A}}^{
*}$ are called  basis and  columns with indices $i_{1} ,\ldots
,i_{r} $ of ${\rm {\bf A}} \in {\bf H}^{m\times n}$ are called
basis  as well.
\end{definition}
The following theorem about the basis rows and  columns from
linear algebra generalize in a straight forward way to
quaternions.

\begin{theorem}  The basis rows of ${\rm {\bf A}}^{
* }{\rm {\bf A}}$ and ${\rm {\bf A}}^{ *}\in {\bf H}^{n\times
m} $ are  left-linearly independent, and the basis columns of
${\rm {\bf A}}^{ *} {\rm {\bf A}}$ and ${\rm {\bf A}} \in {\bf
H}^{m\times n}$ are  right-linearly independent.
\end{theorem}
\begin{theorem}\label{kyrc12} An arbitrary column of ${\rm {\bf A}} \in {\bf H}^{m\times n}$
 is a right linear combination of
its basis columns.
\end{theorem}

{\textit{Proof.}} If columns with indices $i_{1} ,\ldots ,i_{r} $
are the basis columns of ${\rm {\bf A}}$, then the basis principal
minor of ${\rm {\bf A}}^{ *} {\rm {\bf A}}  =: \left( {d_{ij}}
\right)_{n\times n} $ is placed on  crossing of its columns and
rows  with indices $i_{1} ,\ldots ,i_{r} $ as well. Denote by
${\rm {\bf M}}$  the matrix of the basis principal minor.
Supplement it by the ($r + 1$)th row and column consisting of
corresponding entries of the $j$-th row and  column of ${\rm {\bf
A}}^{ *} {\rm {\bf A}}$ respectively. Suppose $j \in {\left\{
{i_{1} ,\ldots ,{\left. {i_{r}} \right\}} } \right.}$. By ${\rm
{\bf D}}_{j} $ denote the obtained matrix.
\[
{\rm {\bf D}}_{j}  =
\begin{pmatrix}
  d_{i_{1} i_{1}} & \cdots & d_{i_{1} i_{r}} &d_{i_{1} j} \\
  \cdots & \cdots & \cdots & \cdots \\
d_{i_{r} i_{1}} & \cdots & d_{i_{r} i_{r}}  & d_{i_{r} j}\\
  d_{j{\kern 1pt} i_{1}}  & \cdots & d_{j{\kern 1pt}
i_{r}}   & d_{j{\kern 1pt} j}
\end{pmatrix}
\]
Since   the Hermitian matrix ${\rm {\bf D}}_{j} $ contains  two
same columns,  we obtain by Corollary \ref{kyrc7} $\det {\rm {\bf
D}}_{j} = {\rm{cdet}} _{j}\, {\rm {\bf D}}_{j} = {\sum\limits_{l =
1}^{r} {L_{i_{l} {\kern 1pt} j} \cdot d_{i_{l} {\kern 1pt} j}} } +
L_{j{\kern 1pt} j} \cdot d_{j{\kern 1pt} j} = 0,$ where $L_{i_{l}
j} $ is the left $i_{l}j$-th cofactor of ${\rm {\bf D}}_{j} $.
Whereas $L_{j{\kern 1pt} j} = \det {\rm {\bf M}} \ne 0$, we get

\begin{equation}
\label{kyr10} d_{j{\kern 1pt} j} = - {\sum\limits_{l = 1}^{r}
{\left( {\det {\rm {\bf M}}} \right)^{ - 1}L_{i_{l} j} \cdot
d_{i_{l} j}} } \quad \forall j \in {\left\{ {i_{1} ,\ldots
,{\left. {i_{r}} \right\}} } \right.}.
\end{equation}

Now suppose that $j \notin {\left\{ {i_{1} ,\ldots ,i_{k} ,i_{k +
1} ,\ldots ,{\left. {i_{r}}  \right\}} } \right.}$ and $i_{k} < j
< i_{k + 1} $. Consider the matrix ${\rm {\bf D}}_{j} $ obtained
from ${\rm {\bf M}}$ by supplementing it by the $j$th row and
 column.
\[
{\rm {\bf D}}_{j} =
\begin{pmatrix}
  d_{i_{1} i_{1}}  & \cdots  & d_{i_{1} i_{k}} & d_{i_{1} j} & d_{i_{1} i_{k + 1}} &\cdots & {d_{i_{1} i_{r}} } \\
  \cdots  & \cdots  & \cdots  & \cdots  & \cdots  & \cdots  & \cdots  \\
  d_{i_{k} i_{1}} &  \cdots & d_{i_{k} i_{k}} & d_{i_{k} j} & d_{i_{k} i_{k + 1}} & \cdots & d_{i_{k} i_{r}} \\
  d_{j{\kern 1pt} i_{1}} & \cdots & d_{j{\kern 1pt}i_{k}} & d_{j{\kern 1pt} j} & d_{j{\kern 1pt} i_{k + 1}
} & \cdots & d_{j{\kern 1pt} i_{r}}\\
  d_{i_{k + 1} i_{1}} & \cdots & d_{i_{k + 1} i_{k}} & d_{i_{k + 1} j} & d_{i_{k + 1} i_{k + 1}} &\cdots & d_{i_{k + 1} i_{r}} \\
 \cdots  & \cdots  & \cdots  & \cdots  & \cdots  & \cdots  & \cdots  \\
 d_{i_{r} i_{1}} &\cdots & d_{i_{r} i_{k}} & d_{i_{r} j} & d_{i_{r} i_{k + 1}} & \cdots & d_{i_{r} i_{r}}
\end{pmatrix}
\]
 The matrix ${\rm {\bf D}}_{j} $ is Hermitian in this case as well. Then we
 have
\[
\det {\rm {\bf D}}_{j} = {\rm{cdet}} _{j} {\rm {\bf D}}_{j} =
{\sum\limits_{l = 1}^{r} {L_{i_{l} {\kern 1pt} j} \cdot d_{i_{l}
{\kern 1pt} j}} }  + L_{j{\kern 1pt} j} \cdot d_{j{\kern 1pt}j}=0.
\]
Since $L_{j{\kern 1pt} j} = \det {\rm {\bf M}}\neq 0$, then
\begin{equation}
\label{kyr11} d_{j{\kern 1pt} j} = - {\sum\limits_{l = 1}^{r}
{\left( {\det {\rm {\bf M}}} \right)^{ - 1}L_{i_{l} j} \cdot
d_{i_{l} j}} }, \quad   \,j \notin {\left\{ {i_{1} ,\ldots ,
i_{r}} \right\} } \subset I_{n}.
\end{equation}
Combining (\ref{kyr10}) and (\ref{kyr11}), we obtain $
 d_{j{\kern 1pt} j} = -
{\sum\limits_{l = 1}^{r} {\left( {\det {\rm {\bf M}}} \right)^{ -
1}L_{i_{l} j} \cdot d_{i_{l} j}} }$,  $ (\forall j = \overline
{1,n}). $ If  $ - \left( {\det {\rm {\bf M}}} \right)^{ -
1}L_{i_{l} j} := \mu _{l}$, then $d_{j{\kern 1pt} j} =
{\sum\limits_{l = 1}^{r} {\mu _{l} \cdot d_{i_{l} j}} }$. Since
$d_{j{\kern 1pt} j} = {\sum\limits_{k = 1}^{m} {\overline
{a_{k{\kern 1pt} j}} a_{k{\kern 1pt} j}} } $ and $d_{i_{l} j} =
{\sum\limits_{k = 1}^{m} {\overline {a_{k{\kern 1pt} {\kern 1pt}
i_{l}} } a_{k{\kern 1pt} j}} }$, then
 \[{\sum\limits_{k = 1}^{m}
{\overline {a_{k{\kern 1pt} j} }  a_{k{\kern 1pt} j}} }  =
{\sum\limits_{l = 1}^{r} {\mu _{l}  {\sum\limits_{k = 1}^{m}
{\overline {a_{k{\kern 1pt} i_{l}} }  a_{k{\kern 1pt} j}} } }}  =
{\sum\limits_{k = 1}^{m} {{\sum\limits_{l = 1}^{r} {\mu _{l}
\overline {a_{k{\kern 1pt} {\kern 1pt} i_{l}} }  a_{k{\kern 1pt}
j}} } }}.\] Hence, $\overline {a_{k{\kern 1pt} j}} =
{\sum\limits_{l = 1}^{r} {\mu _{l} \overline {a_{k{\kern 1pt}
{\kern 1pt} i_{l}} } } } $ and  $a_{k{\kern 1pt} j} =
{\sum\limits_{l = 1}^{r} {a_{k{\kern 1pt} {\kern 1pt} i_{l}}
\overline {\mu _{l}} } }$  $(\forall k = \overline {1,m})$.
Therefore, an arbitrary column of the matrix ${\rm {\bf A}}$ is
the right linear combination of its basis columns with the
coefficients $\overline {\mu _{1}} ,\ldots ,\overline {\mu _{r}}
$, i.e.: ${\rm {\bf a}}_{.\,i_{1}}  \cdot \overline {\mu _{1}} +
\ldots + {\rm {\bf a}}_{.\,i_{r}}  \cdot \overline {\mu _{r}}  =
{\rm {\bf a}}_{.\,j}$\, $(\forall i_{l} \in I_{n}, \, \forall l =
\overline {1,r}).$ $\blacksquare$

The following theorem is proved in a similar manner.
\begin{theorem}\label{kyrc13} An arbitrary row of
${\rm {\bf A}} \in {\bf H}^{m\times n}$ is a left linear
combination of its basis rows.
\end{theorem}
An obvious conclusion of Theorems \ref{kyrc10}, \ref{kyrc11},
\ref{kyrc12}, and \ref{kyrc13} is the criterion of a
nonsingularity of the corresponding Hermitian matrix:
\begin{theorem}\label{kyrc14} The right linearly independence of
columns of ${\rm {\bf A}} \in {\bf H}^{m\times n}$ or the left
linearly independence of rows of ${\rm {\bf A}}^{ *} $ is the
necessary and sufficient condition for $\det {\rm {\bf A}}^{ *}
{\rm {\bf A}} \neq 0.$
\end{theorem}
\section{ Properties of the double determinant \\of a quaternion  square matrix}
Suppose the matrix ${\rm {\bf E}}_{ij}=(e_{pq})_{n\times n}$ such
that $e_{pq}={\left\{ {{\begin{array}{*{20}c}
 {1,\,p=i,\,q=j,} \hfill \\
 {0 ,\,p\neq i, \,q\neq j.} \hfill \\
\end{array}} } \right.}$
\begin{definition}  The matrix ${\rm {\bf P}}_{ij} (b): = {\rm {\bf I}} + b
\cdot {\rm {\bf E}}_{ij}\in {\rm M}(n,{\bf H})$ for $i \ne j$ is
called an elementary unimodular matrix, where ${\rm {\bf I}}$ is
the identity matrix. Matrices ${\rm {\bf P}}_{ij} (b)$ for  $i \ne
j$ and $\forall b \in {\bf H}$ generate the unimodular group ${\rm
{SL}}(n,{\bf H})$, its elements is called the unimodular matrices.
\end{definition}

\begin{theorem}\label{kyrc15} If ${\rm {\bf A}}\in {\rm M}(n,{\bf H})$ is
a Hermitian matrix  and ${\rm {\bf P}}_{ij} \left( {b} \right)$ is
an elementary unimodular matrix, then $\det {\rm {\bf A}} = \det
\left( {{\rm {\bf P}}_{ij} \left( {b} \right) \cdot {\rm {\bf A}}
\cdot {\rm {\bf P}}_{ij}^{ *} \left( {b} \right)} \right).$
\end{theorem}
{\textit{Proof.}} First  note that for $\forall {\rm {\bf U}} \in
{\rm M}(n,{\bf H})$ and a Hermitian matrix ${\rm {\bf A}}$, the
matrix ${\rm {\bf U}}^{ *} {\rm {\bf A}}{\rm {\bf U}}$ is
Hermitian as well. Really, $\left( {{\rm {\bf U}}^{ *} {\rm {\bf
A}}{\rm {\bf U}}} \right)^{ *}  = {\rm {\bf U}}^{ *} {\rm {\bf
A}}^{ *} {\rm {\bf U}} = {\rm {\bf U}}^{ *} {\rm {\bf A}}{\rm {\bf
U}}$. Multiplying a matrix ${\rm {\bf A}}$ by ${\rm {\bf P}}_{ij}
\left( {b} \right)$ on the left adds the $j$th row left-multiplied
by $b$ to the $i$th row. Whereas multiplying a matrix ${\rm {\bf
A}}$ by ${\rm {\bf P}}_{ij}^{ *} \left( {b} \right)$ on the right
adds the $j$th column right-multiplied by $\overline {b}$ to the
$j$th column. Therefore,
\[
{\rm {\bf P}}_{ij} \left( {b} \right) \cdot {\rm {\bf A}} \cdot
{\rm {\bf P}}_{ij}^{ *}  \left( {b} \right) =
\]
\[
  \begin{pmatrix}
  a_{11} & \ldots & a_{1i} + a_{1j} \overline {b} & \ldots & a_{1n} \\
  \ldots &\ldots & \ldots & \ldots & \ldots \\
  a_{i1} + ba_{j1} & \ldots & ( {ba_{jj} + a_{ij}})\overline {b} + ba_{ji} + a_{ii}  & \ldots & a_{in} + ba_{jn} \\
   \ldots &\ldots & \ldots & \ldots & \ldots \\
  a_{n1} & \ldots & a_{ni} + a_{nj}\, \overline {b} & \ldots & a_{nn}
\end{pmatrix}
\]

Then by Theorem \ref{kyrc4} and \ref{kyrc5}, we get
\[
\det \left( {{\rm {\bf P}}_{ij} \left( {b} \right) \cdot {\rm {\bf
A}} \cdot {\rm {\bf P}}_{ij}^{ *}  \left( {b} \right)} \right) =
{\rm{cdet}} _{i} \left( {{\rm {\bf P}}_{ij} \left( {b} \right)
\cdot {\rm {\bf A}} \cdot {\rm {\bf P}}_{ij}^{ *}  \left( {b}
\right)} \right)=\]
\[ ={\rm{cdet}} _{i}\begin{pmatrix}
  a_{11} & ... & a_{1i} & ... & a_{1n} \\
 ... & ... & ... & ... & ...\\
  a_{i1} + ba_{j1} & ... & a_{ii} +  ba_{ji}& ... & a_{in} + ba_{jn} \\
  ... & ... & ... & ... & ...\\
  a_{n1} & ... & a_{ni} & ... & a_{nn}
\end{pmatrix}+\]
\[
+{\rm{cdet}} _{i}\begin{pmatrix}
   a_{11} & ... & a_{1j} \overline {b}  & ... & a_{1n} \\
 ... & ... & ... & ... & ...\\
  a_{i1} + ba_{j1} & ... & ( {ba_{jj} + a_{ij}})\overline {b} & ... & a_{in} + ba_{jn} \\
  ... & ... & ... & ... & ...\\
  a_{n1} & ... & a_{nj}\, \overline {b} & ... & a_{nn}
\end{pmatrix}=\]
\[={\rm{cdet}} _{i} {\rm {\bf A}}+ {\rm{cdet}} _{i} {\rm {\bf
A}}_{i.}  ({b \cdot {\rm {\bf a}}_{j.}})+{\rm{cdet}} _{i} {\rm
{\bf A}}_{.\,i} ({\rm {\bf a}}_{.j})\cdot \overline {b}+\]
\[+{\rm{cdet}} _{i}\begin{pmatrix}
   a_{11} & ... & a_{1j}   & ... & a_{1n} \\
 ... & ... & ... & ... & ...\\
  ba_{j1} & ... &  ba_{jj}  & ... &  ba_{jn} \\
  ... & ... & ... & ... & ...\\
  a_{n1} & ... & a_{nj}& ... & a_{nn}
\end{pmatrix} \cdot \overline {b}.
\]
The matrix
 $
\begin{pmatrix}
   a_{11} & ... & a_{1j}   & ... & a_{1n} \\
 ... & ... & ... & ... & ...\\
  ba_{j1} & ... &  ba_{jj}  & ... &  ba_{jn} \\
  ... & ... & ... & ... & ...\\
  a_{n1} & ... & a_{nj}& ... & a_{nn}
\end{pmatrix}=({\rm {\bf A}}_{.\,i} ({\rm {\bf
a}}_{.j}))_{i\,.} ( b {\rm {\bf a}}_{j\,.})$
 is obtained
from ${\rm {\bf A}}$
 by replacing its $i$th column with the $j$th
column, and then by replacing the $i$th row of the obtained matrix
with its $j$th row left-multiplied by $b$. The $i$th row of (${\rm
{\bf A}}_{.\,i} ({\rm {\bf a}}_{.j}))_{i\,.} ( b {\rm {\bf
a}}_{j\,.})$ is $b {\rm {\bf a}}_{j\,.}$ and its $j$th row is $
{\rm {\bf a}}_{j.}$. Then by Theorem \ref{kyrch24}, we get
${\rm{cdet}} _{i}({\rm {\bf A}}_{.\,i} ({\rm {\bf
a}}_{.j}))_{i\,.} ( b {\rm {\bf a}}_{j\,.})=0$. Furthermore by
Theorem \ref{kyrc8} we have  ${\rm{cdet}} _{i} {\rm {\bf
A}}_{.\,i} ({\rm {\bf a}}_{.j}) =0$, and by Theorem \ref{kyrch24}
we obtain ${\rm{cdet}} _{i} {\rm {\bf A}}_{i\,.} ( b \cdot {\rm
{\bf a}}_{j.})=0.$

Finally, we have $\det \left( {{\rm {\bf P}}_{ij} \left( {b}
\right) \cdot {\rm {\bf A}} \cdot {\rm {\bf P}}_{ij}^{ *}  \left(
{b} \right)} \right)={\rm{cdet}} _{i}{\rm {\bf A}} =\det{\rm {\bf
A}}.$ $\blacksquare$

\begin{theorem}\label{kyrc16} If  ${\rm {\bf A}}\in
 {\rm M}\left( {n,\bf H} \right)$ is a Hermitian matrix
 and $\forall{\rm {\bf U}} \in {\rm SL}(n,{\bf H})$, then
\[
\det {\rm {\bf A}} = \det \left( {{\rm {\bf U}} \cdot {\rm {\bf
A}} \cdot {\rm {\bf U}}^{ *} } \right).
\]
\end{theorem}
{\textit{Proof.}} We claim that  there exists  $ {\left\{ {{\rm
{\bf P}}_{1} ,\ldots ,{\rm {\bf P}}_{k}} \right\}} \subset {\rm
SL}(n,{\bf H}) $ and $ \exists k \in N$ for ${\rm {\bf U}} \in
{\rm SL}\left( {n,\bf H} \right)$ such that ${\rm {\bf U}} = {\rm
{\bf P}}_{k} \cdot \ldots \cdot {\rm {\bf P}}_{1} $. Then ${\rm
{\bf U}}^{ *}  = {\rm {\bf P}}_{1}^{ *}  \cdot \ldots \cdot {\rm
{\bf P}}_{k}^{ *} $.

We prove the theorem by induction on $k$.

 i) The case $k = 1$ has been proved Theorem
\ref{kyrc15}.

 ii) Suppose the  theorem is valid for $k - 1$. That is
${\rm {\bf U}} = {\rm {\bf P}}_{k - 1} \cdot \ldots \cdot {\rm
{\bf P}}_{1} $ and
\[
\det {\rm {\bf A}} = \det \left( {{\rm {\bf P}}_{k - 1} \cdot
\ldots \cdot {\rm {\bf P}}_{1} \cdot {\rm {\bf A}} \cdot {\rm {\bf
P}}_{1} ^{ *}  \cdot \ldots \cdot {\rm {\bf P}}_{k - 1} ^{ *} }
\right){\rm .}
\]
Denote ${\rm {\bf \tilde {A}}}: = {\rm {\bf P}}_{k - 1} \cdot
\ldots \cdot {\rm {\bf P}}_{1} \cdot {\rm {\bf A}} \cdot {\rm {\bf
P}}_{1} ^{ *}  \cdot \ldots \cdot {\rm {\bf P}}_{k - 1} ^{ *} $.
As shown in Theorem \ref{kyrc15} a matrix ${\rm {\bf \tilde {A}}}$
is Hermitian.

iii) If now ${\rm {\bf U}} = {\rm {\bf P}}_{k} \cdot {\rm {\bf
P}}_{k - 1} \ldots \cdot {\rm {\bf P}}_{1} $, then
\[
\det \left( {{\rm {\bf U}} \cdot {\rm {\bf A}} \cdot {\rm {\bf
U}}^{ *} } \right) = \det \left( {{\rm {\bf P}}_{k} \cdot {\rm
{\bf \tilde {A}}} \cdot {\rm {\bf P}}_{k}^{ *} }  \right) = \det
{\rm {\bf \tilde {A}}} = \det {\rm {\bf A}}.\blacksquare
\]

\begin{theorem}\label{kyrc17} If ${\rm {\bf A}}\in {\rm M}\left( {n,\bf H} \right)$
 is a Hermitian matrix, then
 $\exists{\rm {\bf U}} \in {\rm SL}(n,{\bf H})$ and $\exists\mu _{i} \in {\bf R}$, $(\forall i =
\overline {1,n}),$ such that
 ${\rm {\bf U}} \cdot {\rm {\bf A}} \cdot {\rm {\bf
U}}^{ *}  = {\rm {\rm{diag}}}( \mu _{1} ,\ldots ,\mu _{n})$, where
${\rm {\rm{diag}}}( \mu _{1} ,\ldots ,\mu _{n})$ is a diagonal
matrix. Then $\det {\rm {\bf A}} = \mu _{1} \cdot \ldots \cdot \mu
_{n} .$
\end{theorem}
 {\textit{Proof.}} Consider the first column of a Hermitian
 matrix ${\rm {\bf A}}\in {\rm M}\left( {n,\bf H} \right)$. It is
 possible the
following cases.

 i) If $a_{11} \ne 0$, then $\mu
_{1} = a_{11} \in{\bf R}$. By sequentially left-multiplying the
matrix ${\rm {\bf A}}$ by elementary unimodular matrices ${\rm
{\bf P}}_{i1} \left( { - {\frac{{a_{i1}} }{{\mu _{1}} }}}
\right)$, $ \left( {\forall i = \overline {2,n}} \right)$, we
obtain zero for all entries of the first column save for diagonal.
Since $\overline { - {\frac{{a_{i1}} }{{\mu _{1}} }}} = -
{\frac{{a_{1i}} }{{\mu _{1}} }}$, then ${\rm {\bf P}}_{i1}^{ *}
\left( { - {\frac{{a_{i1}} }{{\mu _{1}} }}} \right) = {\rm {\bf
P}}_{1i} \left( { - {\frac{{a_{1i}} }{{\mu _{1}} }}} \right)$. By
sequentially  right-multiplying   the matrix ${\rm {\bf A}}$ by
elementary unimodular matrices ${\rm {\bf P}}_{i1}^{ *} \left( { -
{\frac{{a_{i1}} }{{\mu _{1}} }}} \right)$,
 we get zero for all entries of the first
row save for diagonal. Due to Theorem \ref{kyrc15} the obtained
matrix is Hermitian as well.

ii) Suppose $a_{11} = 0$ and $\exists i \in I_{n}$ $a_{i1} \ne 0$.
Having multiplied  the matrix ${\rm {\bf A}}$
 by elementary unimodular matrices ${\rm {\bf
P}}_{1i} \left( {a_{1i}}  \right)$ on the left and  by ${\rm {\bf
P}}_{i1} \left( {a_{i1}}  \right)$ on the right,  we  get the
matrix ${\rm {\bf \tilde {A}}}$  with an entry $\tilde {a}_{11} =
{\rm n}(a_{i1} )\left( {2 + a_{i\,i}} \right) \in{\bf R}$. Let now
$\mu _{1} = \tilde {a}_{11} $. Again by sequentially multiplying
the matrix ${\rm {\bf \tilde {A}}}$ by ${\rm {\bf P}}_{i1} \left(
{ - {\frac{{\tilde {a}_{i1}} }{{\mu _{1}} }}} \right)$  on the
left and  by ${\rm {\bf P}}_{i1}^{ *} \left( { - {\frac{{ \tilde
{a}_{i1} }}{{\mu _{1}} }}} \right)$, $\left( {\forall i =
\overline {2,n}} \right)$, on the right, we  obtain the matrix
with zero for all entries of the first row and column save for
diagonal.

 iii) If $\forall i \in I_{n}$ $a_{i1} =0$, then put $\mu _{1} =
 a_{11}$.

 Having carried through the described procedure for
all diagonal entries and entries of corresponding rows and columns
by means of a finite  number of multiplications the Hermitian
matrix ${\rm {\bf A}}$   by elementary unimodular matrices ${\rm
{\bf P}}_{k} = {\rm {\bf P}}_{ij} \left( {b_{k}} \right)$ on the
left and  by ${\rm {\bf P}}_{k}^{ *}  = {\rm {\bf P}}_{ji} \left(
{\overline {b_{k}} } \right)$ on the right, we obtain the diagonal
matrix with diagonal entries $\mu _{i} \in {\bf R} \left( {\forall
i = \overline {1,n}} \right)$. Suppose ${\rm {\bf U}} =
{\prod\limits_{k} {{\rm {\bf P}}_{k}} } $, then by Theorem
\ref{kyrc16} we finally obtain
\[
\det ({\rm {\bf U}} \cdot {\rm {\bf A}} \cdot {\rm {\bf U}}^{ *} )
= \det \,(\rm{{\rm{diag}}}\left( {\mu _{1} ,\ldots ,\mu _{n}}
\right)) = \mu _{1} \cdot \ldots \cdot \mu _{n}.\blacksquare
\]
\begin{theorem}\label{kyrc18} If ${\rm {\bf A}} \in {\rm M}\left(
{n,{\bf H}} \right)$, then $\det {\rm {\bf A}}{\rm {\bf A}}^{ *} =
\det {\rm {\bf A}}^{ * }{\rm {\bf A}}$.
\end{theorem}
{\textit{Proof.}} Suppose ${\rm {\bf A}} \in {\rm M}\left( {n,{\bf
H}} \right)$. The matrices $\begin{pmatrix}
  - {\rm {\bf I}} & {\rm {\bf A}} \\
  {\rm {\bf A}}^{ *} & {\rm {\bf 0}}
\end{pmatrix}$ and $\begin{pmatrix}
  {\rm {\bf 0}} & {\rm {\bf A}} \\
  {\rm {\bf A}}^{ *} & - {\rm {\bf I}}
\end{pmatrix}$ are Hermitian. It is easy to see that
\[ \det
\begin{pmatrix}
  - {\rm {\bf I}} & {\rm {\bf A}} \\
  {\rm {\bf A}}^{ *} & {\rm {\bf 0}}
\end{pmatrix}=\det\begin{pmatrix}
  {\rm {\bf 0}} & {\rm {\bf A}} \\
  {\rm {\bf A}}^{ *} & - {\rm {\bf I}}
\end{pmatrix}.\]
The matrix $\begin{pmatrix}
 {\rm {\bf I}} & {\rm {\bf A}} \\
  {\rm {\bf 0}} & {\rm {\bf I}}
\end{pmatrix}$
 can be represented as a product of $n^{2}$ elementary unimodular
$2n\times 2n$ matrices, i.e. $\forall k = \overline {1,n^{2}}$\,
$\exists i = \overline {1,n}  \quad \exists j = \overline {n +
1,n^{2}} $\, $\exists {\rm {\bf P}}_{k} = {\rm {\bf
P}}_{ij}^{\left( {k} \right)} \left( {a_{ij}}  \right)$: $
\begin{pmatrix}
 {\rm {\bf I}} & {\rm {\bf A}} \\
  {\rm {\bf 0}} & {\rm {\bf I}}
\end{pmatrix} = {\prod\limits_{k} {{\rm {\bf P}}_{k}} } $.
Therefore,
$\begin{pmatrix}
 {\rm {\bf I}} & {\rm {\bf A}} \\
  {\rm {\bf 0}} & {\rm {\bf I}}
\end{pmatrix}
 \in {\rm SL}(2n,{\bf H})$. In a similar manner, $\begin{pmatrix}
 {\rm {\bf I}} & {\rm {\bf 0}} \\
  {\rm {\bf A}}^{ *} & {\rm {\bf I}}
\end{pmatrix}\in {\rm SL}(2n,{\bf H})$. From this  by Theorem
\ref{kyrc16}, we have
\[
\begin{array}{c}
 ( - 1)^{n}\det {\rm {\bf A}}{\rm {\bf A}}^{ *}=
   \det\begin{pmatrix}
 {\rm {\bf A}}{\rm {\bf A}}^{ *} &  {\rm {\bf 0}}  \\
  {\rm {\bf 0}} & -{\rm {\bf I}}
\end{pmatrix}=\\
   =\det \left( {\begin{pmatrix}
 {\rm {\bf I}} & {\rm {\bf A}} \\
  {\rm {\bf 0}} & {\rm {\bf I}}
\end{pmatrix}\begin{pmatrix}
  {\rm {\bf 0}} & {\rm {\bf A}} \\
  {\rm {\bf A}}^{ *} & - {\rm {\bf I}}
\end{pmatrix}\begin{pmatrix}
 {\rm {\bf I}} & {\rm {\bf 0}} \\
  {\rm {\bf A}}^{ *} & {\rm {\bf I}}
\end{pmatrix} }\right)=\det\begin{pmatrix}
   {\rm {\bf 0}} & {\rm {\bf A}} \\
  {\rm {\bf A}}^{ *} & - {\rm {\bf I}}
\end{pmatrix}=\\
=\det\begin{pmatrix}
  - {\rm {\bf I}} & {\rm {\bf A}} \\
  {\rm {\bf A}}^{ *} & {\rm {\bf 0}}
\end{pmatrix}
   =\det\left( {\begin{pmatrix}
 {\rm {\bf I}} & {\rm {\bf 0}} \\
  {\rm {\bf A}}^{ *} & {\rm {\bf I}}
\end{pmatrix}\begin{pmatrix}
  - {\rm {\bf I}} & {\rm {\bf A}} \\
  {\rm {\bf A}}^{ *} & {\rm {\bf 0}}
\end{pmatrix}\begin{pmatrix}
 {\rm {\bf I}} & {\rm {\bf A}} \\
  {\rm {\bf 0}} & {\rm {\bf I}}
\end{pmatrix}}\right)=\\
  =\det\begin{pmatrix}
 {\rm {\bf A}}{\rm {\bf A}}^{ *} &  {\rm {\bf 0}}\\
  {\rm {\bf 0}} & -{\rm {\bf I}}
\end{pmatrix}=
( - 1)^{n}\det {\rm {\bf A}}^{ *}{\rm {\bf A}}. \blacksquare
\end{array}
\]
\begin{definition} For $\forall {\rm {\bf A}} \in {\rm M}\left(
{n,{\bf H}} \right)$ the determinant of its corresponding
Hermitian matrix is called its double determinant, i.e.
\[{\rm{ddet}}{ \rm{\bf A}}: = \det \left( {{\rm {\bf A}}^{ *} {\rm
{\bf A}}} \right) = \det \left( {{\rm {\bf A}}{\rm {\bf A}}^{ *} }
\right).\]
\end{definition}
\begin{theorem} \label{kyrc19} If $\forall{\left\{
{{\rm {\bf A}},{\rm {\bf B}}} \right\}} \subset {\rm M}\left(
{n,{\bf H}} \right)$, then ${\rm{ddet}} \left( {{\rm {\bf A}}
\cdot {\rm {\bf B}}} \right) = {\rm{ddet}}{ \rm{\bf A}}
\cdot{\rm{ddet}} {\rm {\bf B}}.$
\end{theorem}
{\textit{Proof.}} Due to Theorem \ref{kyrc17} for the Hermitian
matrix ${\rm {\bf A}}^{ *} {\rm {\bf A}}$, there exists ${\rm {\bf
U}} \in {\rm SL}\left( {n,{\bf H}} \right)$ such that ${\rm {\bf
U}}^{ *} \cdot {\rm {\bf A}}^{ *} {\rm {\bf A}} \cdot {\rm {\bf
U}}= \left( {{\rm {\bf A}} \cdot {\rm {\bf U}}} \right)^{ *} \cdot
{\rm {\bf A}} \cdot {\rm {\bf U}}={\rm{{\rm{diag}}}} \left(
{\alpha _{1} ,\ldots ,\alpha _{n}} \right)$, where $\alpha _{i}
\in{\bf R}$. If ${\rm {\bf A}} \cdot {\rm {\bf U}}=(q_{i\,j})_{n
\times n}$, then $\alpha _{i} = {\sum\limits_{k} {\overline
{q_{ki}}  q_{ki} = {\sum\limits_{k} {n(q_{ki} )}} }} \in {\bf R}_{
+}$\, $ ({\forall i = \overline {1,n}})$, where ${\bf R}_{ +}$ is
the set of the nonnegative real numbers.
 Therefore for $\forall \alpha _{i}\in {\bf R}_{ +}\,\,  \exists
\sqrt {\alpha _{i} } \in {\bf R}_{ +} $\, $ \left( {\forall i =
\overline {1,n}} \right)$. By virtue of $({\rm {\bf U}}^{ *} )^{ -
1} = ({\rm {\bf U}}^{ - 1})^{ *} $ for Hermitian $({\rm {\bf U}}^{
- 1}{\rm {\bf B}})^{ *} ({\rm {\bf U}}^{ - 1}{\rm {\bf B}})$ there
exist ${\rm {\bf W}} \in {\rm SL}\left( {n,{\bf H}} \right)$ and $
\beta _{i} \in{\bf R_{ +}}$ $ \left( {\forall i = \overline {1,n}}
\right)$ such that $
 {\rm {\bf W}}^{ *} ({\rm {\bf U}}^{ - 1}{\rm {\bf
B}})^{ *} ({\rm {\bf U}}^{ - 1}{\rm {\bf B}}){\rm {\bf W}} =
{\rm{diag}}(\beta _{1} ,\ldots ,\beta _{n} )$.  Hence by Theorems
\ref{kyrc17} and \ref{kyrc18}, we obtain
\[
\begin{array}{c}
 {\rm{ddet}} ({\rm {\bf A}} \cdot {\rm {\bf B}}) = \det ({\rm {\bf
B}}^{ *} ({\rm {\bf A}}^{ *} {\rm {\bf A}}) {\rm {\bf B}}) = \det
({\rm {\bf B}}^{ *} ({\rm {\bf U}}^{ *} )^{ - 1}{\rm {\bf U}}^{ *}
({\rm {\bf A}}^{ *} {\rm {\bf A}}){\rm {\bf U}}{\rm {\bf U}}^{ -
1}{\rm {\bf B}}) = \\
   = \det \left( {\left( {{\rm {\bf U}}^{ - 1}{\rm {\bf B}}} \right)^{ *
}{\rm{diag}}\left( {\alpha _{1} ,\ldots ,\alpha _{n}}  \right){\rm
{\bf U}}^{ - 1}{\rm {\bf B}}} \right) = \\
    = \det \left( {({\rm{diag}}\left( {\sqrt {\alpha _{1}}  ,\ldots ,\sqrt {\alpha
_{n}} }  \right){\rm {\bf U}}^{ - 1}{\rm {\bf B}})^{ *}
({\rm{diag}}\left( {\sqrt {\alpha _{1}}  ,\ldots ,\sqrt {\alpha
_{n}} }  \right){\rm {\bf U}}^{ - 1}{\rm {\bf B}})} \right) =\\
   = \det \left( {({\rm{diag}}\left( {\sqrt {\alpha _{1}}  ,\ldots ,\sqrt {\alpha
_{n}} }  \right){\rm {\bf U}}^{ - 1}{\rm {\bf
B}})({\rm{diag}}\left( {\sqrt {\alpha _{1}} ,\ldots ,\sqrt {\alpha
_{n}} }  \right){\rm {\bf U}}^{ - 1}{\rm {\bf B}})^{ *} } \right)
=\\
    = \det \left( {{\rm{diag}}\left( {\sqrt {\alpha _{1}}  ,\ldots ,\sqrt {\alpha _{n}
}}  \right)({\rm {\bf U}}^{ - 1}{\rm {\bf B}})({\rm {\bf U}}^{ -
1}{\rm {\bf B}})^{ *} {\rm{diag}}\left( {\sqrt {\alpha _{1}}
,\ldots ,\sqrt {\alpha _{n}} } \right)} \right) =\\
 = \det \left(
{{\rm{diag}}\left( {\sqrt {\alpha _{1}}  ,\ldots ,\sqrt {\alpha
_{n} }}  \right) ({\rm {\bf W}}^{ - 1})^{ *} {\rm{diag}}\left(
{\beta _{1} ,\ldots ,\beta _{n}} \right) {\rm {\bf W}}^{ - 1}
\times }\right.\\
 \left.{ \times {\rm{diag}}\left( {\sqrt {\alpha
_{1}} ,\ldots ,\sqrt {\alpha _{n}} }  \right) }\right) =
 \det \left( {   \left( {({\rm {\bf W}}^{ - 1})^{T}}
\right)^{ *}   {\rm{diag}}\left( {\sqrt {\alpha _{1}} ,\ldots
,\sqrt {\alpha _{n}} }\right) }\right.\times\\
 \left.{
\times {\rm{diag}}\left( {\beta _{1} ,\ldots ,\beta _{n}}  \right)
\cdot{\rm{diag}}\left( {\sqrt {\alpha _{1}}  ,\ldots ,\sqrt
{\alpha _{n}} } \right) ({\rm {\bf W}}^{ - 1})^{T}   } \right) =\\
    \det \left( {{\rm{diag}}\left( {\sqrt {\alpha _{1}}  ,\ldots ,\sqrt {\alpha _{n}
}}  \right) \cdot{\rm{diag}}\left( {\beta _{1} ,\ldots ,\beta
_{n}} \right){\rm{diag}}\left( {\sqrt {\alpha _{1}} ,\ldots ,\sqrt
{\alpha _{n}} }  \right)} \right) =\\
    = \alpha _{1} \cdot \ldots \cdot \alpha _{n} \cdot \beta _{1} \cdot \ldots
\cdot \beta _{n}  = \det {\rm {\bf A}} \cdot \det {\rm {\bf B}} =
\det {\rm {\bf B}} \cdot \det {\rm {\bf A}}. \blacksquare\\
\end{array}
\]
\begin{remark} The proofs of Theorems  \ref{kyrc18} and \ref{kyrc19} are similarly
 to the proofs in \cite[p.533]{ch3},
 and they differ  by using  different determinant
functionals.
\end{remark}
\begin{remark} From Theorems \ref{kyrc14} and \ref{kyrc19} follows that
 for  $\forall {\rm {\bf A}} \in
{\rm M}\left( {n,{\bf H}} \right)$ ${\rm{ddet}} {\rm {\bf A}}$
satisfies
 Axioms 1, 2, 3. From {\cite{As1,co4}} we have
\[{\rm{ddet}} {\rm {\bf A}} = {\rm{Mdet}}\, \left( {{\rm {\bf
A}}^{
* }{\rm {\bf A}}} \right) =  {\rm{Sdet}} {\rm {\bf A}} =
{\rm{Ddet}} ^{2 }{\rm {\bf A}},\] where ${\rm{Sdet}} {\rm {\bf
A}}$, ${\rm{Ddet}} {\rm {\bf A}}$ are accordingly the determinants
of Study and  of Diedonne.
\end{remark}

\section{Determinant representation  \\ of the quaternion inverse matrix}
\begin{definition} Suppose a matrix ${\rm {\bf A}} \in {\rm M}(n,{\bf
H})$ and ${\rm{ddet}} {\rm {\bf A}} = {\rm{cdet}} _{j} \left(
{{\rm {\bf A}}^{ *} {\rm {\bf A}}} \right) = {\sum\limits_{i}
{{\mathbb{L}} _{ij} \cdot a_{ij}} }$, $ (\forall j = \overline
{1,n} )$, then ${\mathbb{L}} _{ij} $ is called the left double
$ij$-th cofactor of ${\rm {\bf A}}$.
\end{definition}
\begin{definition} Suppose a matrix ${\rm {\bf A}} \in {\rm M}(n,{\bf
H})$ and ${\rm{ddet}} {\rm {\bf A}} = {\rm{rdet}}_{i} \left( {{\rm
{\bf A}}{\rm {\bf A}}^{ *} } \right) = {\sum\limits_{j} {a_{ij}
\cdot} } {\mathbb{R}} _{ i{\kern 1pt}j}$, $ (\forall i = \overline
{1,n}) $, then ${\mathbb{R}} _{ i{\kern1pt}j} $ is called
 the right double $ij$th cofactor  of ${\rm {\bf A}}$.
\end{definition}
\begin{theorem} \label{kyrc20} The necessary and sufficient condition of invertibility
of  ${\rm {\bf A}} \in {\rm M}(n,{\bf H})$ is ${\rm{ddet}} {\rm
{\bf A}} \ne 0$. Then $\exists {\rm {\bf A}}^{ - 1} = \left(
{L{\rm {\bf A}}} \right)^{ - 1} = \left( {R{\rm {\bf A}}}
\right)^{ - 1}$, where
\begin{equation}
\label{kyr12} \left( {L{\rm {\bf A}}} \right)^{ - 1} =\left( {{\rm
{\bf A}}^{ *}{\rm {\bf A}} } \right)^{ - 1}{\rm {\bf A}}^{ *}
={\frac{{1}}{{{\rm{ddet}}{ \rm{\bf A}} }}}
\begin{pmatrix}
  {\mathbb{L}} _{11} & {\mathbb{L}} _{21}& \ldots & {\mathbb{L}} _{n1} \\
  {\mathbb{L}} _{12} & {\mathbb{L}} _{22} & \ldots & {\mathbb{L}} _{n2} \\
  \ldots & \ldots & \ldots & \ldots \\
 {\mathbb{L}} _{1n} & {\mathbb{L}} _{2n} & \ldots & {\mathbb{L}} _{nn}
\end{pmatrix}
\end{equation}
\begin{equation}\label{kyr13} \left( {R{\rm {\bf A}}} \right)^{ - 1} = {\rm {\bf
A}}^{ *} \left( {{\rm {\bf A}}{\rm {\bf A}}^{ *} } \right)^{ - 1}
= {\frac{{1}}{{{\rm{ddet}}{ \rm{\bf A}}^{ *} }}}
\begin{pmatrix}
 {\mathbb{R}} _{\,{\kern 1pt} 11} & {\mathbb{R}} _{\,{\kern 1pt} 21} &\ldots & {\mathbb{R}} _{\,{\kern 1pt} n1} \\
 {\mathbb{R}} _{\,{\kern 1pt} 12} & {\mathbb{R}} _{\,{\kern 1pt} 22} &\ldots & {\mathbb{R}} _{\,{\kern 1pt} n2}  \\
 \ldots  & \ldots & \ldots & \ldots \\
 {\mathbb{R}} _{\,{\kern 1pt} 1n} & {\mathbb{R}} _{\,{\kern 1pt} 2n} &\ldots & {\mathbb{R}} _{\,{\kern 1pt} nn}
\end{pmatrix}
\end{equation}
and ${\mathbb{L}} _{ij} = {\rm{cdet}} _{j} ({\rm {\bf
A}}^{\ast}{\rm {\bf A}})_{.j} \left( {{\rm {\bf a}}_{.{\kern 1pt}
i}^{ *} } \right)$, ${\mathbb{R}} _{\,{\kern 1pt} ij} =
{\rm{rdet}}_{i} ({\rm {\bf A}}{\rm {\bf A}}^{\ast})_{i.} \left(
{{\rm {\bf a}}_{j.}^{ *} }  \right), \left( {\forall i,j =
\overline {1,n}} \right).$ \end{theorem}
 {\textit{Proof.}}
(Necessity). Suppose there exists the inverse matrix  ${\rm {\bf
A}}^{ - 1}$ of
 ${\rm {\bf A}} \in {\rm M}(n,{\bf H})$. By virtue of
$ {\rm rank}\,{\rm {\bf A}} \ge {\rm rank}({\rm {\bf A}}^{ -
1}{\rm {\bf A}}) = {\rm rank}\,{\rm {\bf I}} = n, $ then ${\rm
rank}\,{\rm {\bf A}} = n$. Thus the columns of  $ {\rm {\bf A}}$
are right linearly independent. By Theorem \ref{kyrc14}, this
implies  $\det {\rm {\bf A}}^{ *} {\rm {\bf A}} = {\rm ddet} {\rm
{\bf A}} \ne 0$.

(Sufficiency) Since ${\rm{ddet}}{ \rm{\bf A}} = \det {\rm {\bf
A}}^{ *} {\rm {\bf A}} \ne 0$,  by Theorem \ref{kyrc9}  there
exists the inverse $\left( {{\rm {\bf A}}^{ *} {\rm {\bf A}}}
\right)^{ - 1}$ of the Hermitian matrix ${\rm {\bf A}}^{ *} {\rm
{\bf A}}$. Multiplying it on the right  by ${\rm {\bf A}}^{ *} $
obtains the left  inverse $\left( {L{\rm {\bf A}}} \right)^{ - 1}
= \left( {{\rm {\bf A}}^{ *} {\rm {\bf A}}} \right)^{ - 1}{\rm
{\bf A}}^{ *} $. By representing  $\left( {{\rm {\bf A}}^{ *} {\rm
{\bf A}}} \right)^{ - 1}=({\frac{{L_{ij}}}{{{\rm{ddet}} {\rm {\bf
A}}}}})_{n\times n}$ as the left inverse matrix, we get
\[
   \left( {L{\rm {\bf A}}} \right)^{ - 1} = \left( {L\left( {{\rm {\bf A}}^{ *
}{\rm {\bf A}}} \right)} \right)^{ - 1}{\rm {\bf A}}^{ *}=\]
  \[
 = {\frac{{1}}{{{\rm{ddet}} {\rm {\bf A}}}}}\left( {{\begin{array}{*{20}c}
 {{\sum\limits_{k} {L_{k1} a_{k1}^{ *} } } } \hfill & {{\sum\limits_{k}
{L_{k1} a_{k2}^{ *} } } } \hfill & {\ldots}  \hfill &
{{\sum\limits_{k} {L_{k1} a_{kn}^{ *} } } } \hfill \\
 {{\sum\limits_{k} {L_{k2} a_{k1}^{ *} } } } \hfill & {{\sum\limits_{k}
{L_{k2} a_{k2}^{ *} } } } \hfill & {\ldots}  \hfill &
{{\sum\limits_{k} {L_{k2} a_{kn}^{ *} } } } \hfill \\
 {\ldots}  \hfill & {\ldots}  \hfill & {\ldots}  \hfill & {\ldots}  \hfill
\\
 {{\sum\limits_{k} {L_{kn} a_{k1}^{ *} } } } \hfill & {{\sum\limits_{k}
{L_{kn} a_{k2}^{ *} } } } \hfill & {\ldots}  \hfill &
{{\sum\limits_{k} {L_{kn} a_{kn}^{ *} } } } \hfill \\
\end{array}} } \right) =
\]
\[ = {\frac{{1}}{{{\rm{ddet}}{ \rm{\bf A}}}}}
   \begin{pmatrix}
   {{\rm{cdet}} _{1} ({\rm {\bf
A}}^{\ast}{\rm {\bf A}})_{.\,1} \left( {{\rm {\bf a}}_{.\,1}^{ *}
}  \right)}
 & {{\rm{cdet}} _{1} ({\rm {\bf
A}}^{\ast}{\rm {\bf A}})_{.\,1} \left( {{\rm {\bf a}}_{.\,2}^{ *}
} \right)}  & {\ldots}   & {{\rm{cdet}} _{1} ({\rm {\bf
A}}^{\ast}{\rm {\bf A}})_{.\,1} \left( {{\rm {\bf a}}_{.\,n}^{ *}
} \right)}  \\
   {{\rm{cdet}} _{2} ({\rm {\bf
A}}^{\ast}{\rm {\bf A}})_{.\,2} \left( {{\rm {\bf a}}_{.\,1}^{ *}
} \right)}
 & {{\rm{cdet}} _{2} ({\rm {\bf
A}}^{\ast}{\rm {\bf A}})_{.\,2} \left( {{\rm {\bf a}}_{.\,2}^{ *}
} \right)}  & \ldots & {{\rm{cdet}} _{2} ({\rm {\bf A}}^{\ast}{\rm
{\bf A}})_{.\,2} \left( {{\rm {\bf a}}_{.\,n}^{ *}
} \right)}  \\
\ldots & \ldots & \ldots & \ldots \\
{{\rm{cdet}} _{n} ({\rm {\bf A}}^{\ast}{\rm {\bf A}})_{.\,n}
\left( {{\rm {\bf a}}_{.\,1}^{ *} }  \right)}
 & {{\rm{cdet}} _{n} ({\rm {\bf
A}}^{\ast}{\rm {\bf A}})_{.\,n} \left( {{\rm {\bf a}}_{.\,2}^{ *}
} \right)}  & \ldots & {{\rm{cdet}} _{n} ({\rm {\bf A}}^{\ast}{\rm
{\bf A}})_{.\,n} \left( {{\rm {\bf a}}_{.\,n}^{ *} } \right)}
\end{pmatrix}
\]
By virtue of ${\rm{ddet}}{ \rm{\bf A}} = \det ({\rm {\bf
A}}^{\ast}{\rm {\bf A}}) = {\rm{cdet}} _{j} ({\rm {\bf
A}}^{\ast}{\rm {\bf A}}) = {\sum\limits_{i} {{\rm{cdet}} _{j}
({\rm {\bf A}}^{\ast}{\rm {\bf A}})_{.j} \left( {{\rm {\bf
a}}_{.i}^{ *} } \right)}}  \cdot a_{ij} = {\sum\limits_{i}
{{\mathbb{L}} _{ij} \cdot a_{ij}} } $,  $ \left( {\forall j =
\overline {1,n}} \right)$, we obtain (\ref{kyr12}).

Now prove the formula (\ref{kyr13}). By Theorem \ref{kyrc9} there
exists an inverse matrix $\left( {{\rm {\bf A}} {\rm {\bf A}}^{
*}} \right)^{ - 1}=({\frac{{R_{ij}}}{{{\rm{ddet}} {\rm {\bf
A}}}}})_{n\times n}$. By having left-multiplied it by ${\rm {\bf
A}}^{ *} $, we obtain:
\[
   \left( {R{\rm {\bf A}}} \right)^{ - 1} = {\rm {\bf A}}^{ *} \left(
{R\left( {{\rm {\bf A}}^{ *} {\rm {\bf A}}} \right)} \right)^{ -
1}=\]
\[
 = {\frac{{1}}{{d\det {\rm {\bf A}}}}}\left( {{\begin{array}{*{20}c}
 {{\sum\limits_{k} {a_{1k}^{ *}  R_{1k}} } } \hfill & {{\sum\limits_{k}
{a_{1k}^{ *}  R_{2k}} } } \hfill & {\ldots}  \hfill &
{{\sum\limits_{k} {a_{1k}^{ *}  R_{nk}} } } \hfill \\
 {{\sum\limits_{k} {a_{2k}^{ *}  R_{1k}} } } \hfill & {{\sum\limits_{k}
{a_{2k}^{ *}  R_{2k}} } } \hfill & {\ldots}  \hfill &
{{\sum\limits_{k} {a_{2k}^{ *}  R_{nk}} } } \hfill \\
 {\ldots}  \hfill & {\ldots}  \hfill & {\ldots}  \hfill & {\ldots}  \hfill
\\
 {{\sum\limits_{k} {a_{nk}^{ *}  R_{1k}} } } \hfill & {{\sum\limits_{k}
{a_{1k}^{ *}  R_{1k}} } } \hfill & {\ldots}  \hfill &
{{\sum\limits_{k} {a_{nk}^{ *}  R_{nk}} } } \hfill \\
\end{array}} } \right)
\]
 \[ = {\frac{{1}}{{{\rm{ddet}}{ \rm{\bf A}}^{ *}
}}}\begin{pmatrix}
  {{\rm{rdet}} _{1} ({\rm {\bf A}}{\rm {\bf A}}^{\ast})_{1.}
  \left( {{\rm {\bf a}}_{1.}^{ *} }  \right)}
\hfill & {{\rm{rdet}} _{2} ({\rm {\bf A}}{\rm {\bf
A}}^{\ast})_{2.} \left( {{\rm {\bf a}}_{1.}^{ *} } \right)} \hfill
& {\ldots}  \hfill & {{\rm{rdet}} _{n} ({\rm {\bf A}}{\rm {\bf
A}}^{\ast})_{n.} \left( {{\rm {\bf a}}_{1.}^{ *} } \right)} \hfill
\\
  {{\rm{rdet}} _{1} ({\rm {\bf A}}{\rm {\bf A}}^{\ast})_{1.}
   \left( {{\rm {\bf a}}_{2.}^{ *} }  \right)}
\hfill & {{\rm{rdet}} _{2} ({\rm {\bf A}}{\rm {\bf
A}}^{\ast})_{2.} \left( {{\rm {\bf a}}_{2.}^{ *} } \right)} \hfill
& {\ldots}  \hfill & {{\rm{rdet}} _{n} ({\rm {\bf A}}{\rm {\bf
A}}^{\ast})_{n.} \left( {{\rm {\bf a}}_{2.}^{ *} } \right)} \hfill
\\
    {\ldots}   & {\ldots}  & {\ldots}   & {\ldots}\\
    {{\rm{rdet}} _{1} ({\rm {\bf A}}{\rm {\bf A}}^{\ast})_{1.}
     \left( {{\rm {\bf a}}_{n.}^{ *} }  \right)}
\hfill & {{\rm{rdet}} _{2} ({\rm {\bf A}}{\rm {\bf
A}}^{\ast})_{2.} \left( {{\rm {\bf a}}_{n.}^{ *} } \right)} \hfill
& {\ldots}  \hfill & {{\rm{rdet}} _{n} ({\rm {\bf A}}{\rm {\bf
A}}^{\ast})_{n.} \left( {{\rm {\bf a}}_{n.}^{ *} } \right)} \hfill
  \end{pmatrix}
\]
By virtue of ${\rm{ddet}}{ \rm{\bf A}}  = {\rm{rdet}}_{i} ({\rm
{\bf A}}{\rm {\bf A}}^{\ast}) = {\sum\limits_{j} {a_{ij} \cdot
{\rm{rdet}} _{i} ({\rm {\bf A}}{\rm {\bf A}}^{\ast})_{i.} \left(
{{\rm {\bf a}}_{j.}^{ *} } \right)}}  = {\sum\limits_{j} {a_{ij}
\cdot {\mathbb{R}} _{{\kern 1pt} {\kern 1pt} ij}} }  $, $ \left(
{\forall i = \overline {1,n}} \right)$, the formula (\ref{kyr13})
is valid.  The equality $\left( {L{\rm {\bf A}}} \right)^{ - 1} =
\left( {R{\rm {\bf A}}} \right)^{ - 1}$ is immediately from the
well-known fact that if there exists an inverse matrix over an
arbitrary skew field, then it is unique.$\blacksquare$
\begin{remark} In Theorem \ref{kyrc20} on the
assumption that ${\rm{ddet}}{ \rm{\bf A}} \ne 0$,  the inverse
matrix ${\rm {\bf A}}^{ - 1}$ of  ${\rm {\bf A}} \in {\rm
M}(n,{\bf H})$  is represented by some "double"
 analogue of the classical adjoint matrix. If we denote this
 analogue
 by ${\rm Adj}{\left[ {{\left[ {{\rm
{\bf A}}} \right]}} \right]}$, then the following formula is valid
over
 $\bf H$:
\[
{\rm {\bf A}}^{ - 1} = {\frac{{{\rm Adj}{\left[ {{\left[ {{\rm
{\bf A}}} \right]}} \right]}}}{{{\rm{ddet}} {\rm {\bf A}}}}}.
\]
\end{remark}

\section{ Cramer rule  for quaternionic systems \\ of linear equations}
\begin{theorem} \label{kyrc21}Let
\begin{equation}
\label{kyr14} {\rm {\bf A}} \cdot {\rm {\bf x}} = {\rm {\bf y}}
\end{equation}
be a right system of  linear equations  with a matrix of
coefficients ${\rm {\bf A}}\in {\rm M}(n,{\bf H})$, a column of
constants ${\rm {\bf y}} = \left( {y_{1} ,\ldots ,y_{n} }
\right)^{T}\in {\rm \bf H}^{n\times 1}$, and a column of unknowns
${\rm {\bf x}} = \left( {x_{1} ,\ldots ,x_{n}} \right)^{T}$. If
${\rm{ddet}}{ \rm{\bf A}} \ne 0$, then the solution to the linear
system (\ref{kyr14}) is given by components
\begin{equation}
\label{kyr15} x_{j}^{} = {\frac{{{\rm{cdet}} _{j} ({\rm {\bf A}}^{
*} {\rm {\bf A}})_{.j} \left( {{\rm {\bf f}}}
\right)}}{{{\rm{ddet}} {\rm {\bf A}}}}}, \left( {\forall j =
\overline {1,n}}  \right),
\end{equation}
where $ {\rm {\bf f}} = {\rm {\bf A}}^{ *} {\rm {\bf y}}.$
\end{theorem}
{\textit{Proof.}} By Theorem \ref{kyrc20} ${\rm {\bf A}}$ is
invertibility. So there exists  the unique inverse matrix ${\rm
{\bf A}}^{ - 1}$. From this the existence and  uniqueness of
solutions of (\ref{kyr14}) follows immediately. Consider ${\rm
{\bf A}}^{ - 1}$ as the left inverse $\left( {L{\rm {\bf A}}}
\right)^{ - 1} = \left( {{\rm {\bf A}}^{ *} {\rm {\bf A}}}
\right)^{ - 1}{\rm {\bf A}}^{ *} $. Then  we get $ {\rm {\bf x}} =
{\rm {\bf A}}^{ - 1} \cdot {\rm {\bf y}} = \left( {{\rm {\bf A}}^{
*} {\rm {\bf A}}} \right)^{ - 1}{\rm {\bf A}}^{ *}  \cdot {\rm
{\bf y}}$. Denote ${\rm {\bf f}}: = {\rm {\bf A}}^{ *}  \cdot {\rm
{\bf y}}$. Here ${\rm {\bf f}} = \left( {{\begin{array}{*{20}c}
 {f_{1}}  \hfill & {f_{2}}  \hfill & {\ldots}  \hfill & {f_{n}}  \hfill \\
\end{array}} } \right)^{T}$ is the $n$-dimension column vector over $\bf H$.
By considering  $ \left( {{\rm {\bf A}}^{ *} {\rm {\bf A}}}
\right)^{ - 1}$ as the left inverse, the solution of (\ref{kyr14})
is represented  by components
\[
x_{j} = \left( {{\rm{ddet}}{ \rm{\bf A}}} \right)^{ -
1}{\sum\limits_{i = 1}^{n} {L_{ij} \cdot f_{i}} }\, , \quad \left(
{\forall j = \overline {1,n}}  \right),
\]
where $L_{ij} $ is the left $ij$-th cofactor of the Hermitian
matrix $ \left( {{\rm {\bf A}}^{ *} {\rm {\bf A}}} \right)$. From
here
 we obtain (\ref{kyr15}).$\blacksquare$
\begin{theorem} Let
\begin{equation}\label{kyr16}
 {\rm {\bf x}} \cdot {\rm {\bf A}} = {\rm {\bf y}}
\end{equation}
be a left  system of  linear equations
 with a
matrix of coefficients ${\rm {\bf A}}\in {\rm M}(n,{\bf H})$, a
row of constants ${\rm {\bf y}} = \left( {y_{1} ,\ldots ,y_{n} }
\right)\in {\rm \bf H}^{1\times n}$, and a row of unknowns ${\rm
{\bf x}} = \left( {x_{1} ,\ldots ,x_{n}}  \right)$. If
${\rm{ddet}}{ \rm{\bf A}} \ne 0$, then the solution to the linear
system (\ref{kyr16}) is given by components
\begin{equation}\label{kyr17}
  x_{i}^{} = {\frac{{{\rm{rdet}}_{i} ({\rm {\bf A}}{\rm {\bf A}}^{ *})_{i.} \left( {{\rm
{\bf z}}} \right)}}{{{\rm{ddet}}{ \rm{\bf A}}}}}, \quad \left(
{\forall i = \overline {1,n}}  \right)
\end{equation}
where $ {\rm {\bf z}} = {\rm {\bf y}}{\rm {\bf A}}^{ *} .$
\end{theorem}
{\textit{Proof}}. The proof of this theorem is analogous to that
of Theorem \ref{kyrc21}.
\begin{remark} The formulas (\ref{kyr15}) and (\ref{kyr17}) are the
obvious and natural generalizations of  Cramer rule for
quaternionic systems of  linear equations. The best similarity to
 Cramer rule can be received by Theorem \ref{kyrc9} in the following specific cases.
\end{remark}
\begin{theorem}If the matrix of coefficients  ${\rm {\bf A}}\in {\rm M}(n,{\bf H})$
 in the right system of  linear equations over
$\bf H$  (\ref{kyr14}) is Hermitian, then the unique solution
vector ${\rm {\bf x}}=(x_{1},x_{2},\ldots,x_{n})$  of the system
is given by
\[
x_{j} = {\frac{{{\rm{cdet}} _{j} {\rm {\bf A}}_{.j} \left( {{\rm
{\bf y}}} \right)}}{{\det {\rm {\bf A}}}}} \quad  {\forall j =
\overline {1,n}}.
\]
\end{theorem}
\begin{theorem} If the matrix of coefficients
${\rm {\bf A}}\in {\rm M}(n,{\bf H})$ in the left system of linear
equations over $\bf H$ (\ref{kyr16}) is Hermitian, then the unique
solution vector ${\rm {\bf x}}=(x_{1},x_{2},\ldots,x_{n})$  is
given by
\[
x_{i}^{} = {\frac{{{\rm{rdet}}_{i} {\rm {\bf A}}_{i.} \left( {{\rm
{\bf y}}} \right)}}{{\det {\rm {\bf A}}}}} \quad  {\forall i =
\overline {1,n}}.
\]
\end{theorem}

\end{document}